\documentclass[a5paper,11pt]{amsart}
\pdfoutput=1

\usepackage[finnish, main=english]{babel}
\usepackage[utf8]{inputenc}
\usepackage{bm}
\usepackage{mathtools,amssymb}
\usepackage{hyperref}
\usepackage{xcolor,graphicx}
\usepackage{mathrsfs}
\usepackage[shortlabels]{enumitem}
\usepackage{lineno}
\usepackage{amsmath,eufrak}
\usepackage{enumitem}
\usepackage{amsthm} 
\usepackage{caption}
\usepackage{pdfpages}

\allowdisplaybreaks

\mathtoolsset{showonlyrefs}

\graphicspath{{images/}}

\newtheorem{theorem}{Theorem}[section]

\newcommand{\C}{{\mathbb C}}
\newcommand{\R}{{\mathbb R}}
\newcommand{\Z}{{\mathbb Z}}
\newcommand{\N}{{\mathbb N}}

\newcommand{\der}{{\mathrm d}}
\newcommand*{\sol}[1]{#1^\mathrm{s}}

\newcommand{\riesz}{I_{\alpha}}
\newcommand{\xrt}{X}
\newcommand{\no}{N}

\newcommand{\tempered}{\mathscr{S}^{\prime}}
\newcommand{\rapidly}{\mathscr{O}_C^{\prime}}

\newcommand{\fourier}{\mathcal{F}}
\newcommand{\ifourier}{\mathcal{F}^{-1}}

\newcommand{\cdistr}{\mathcal{E}'}
\newcommand{\distr}{\mathcal{D}^{\prime}}
\newcommand{\dimens}{n}
\newcommand{\kernel}{h_{\alpha}} 
\newcommand{\fraclaplace}{(-\Delta)^s}
\newcommand{\pdos}{\mathcal{P}}
\newcommand{\adm}{\mathcal{A}}

\newcommand{\dplane}{\mathcal{R}_d}
\newcommand{\nod}{\mathcal{N}_d}
\newcommand{\geod}{\mathcal{I}}
\newcommand{\mixed}{L_{k, l}}

\newcommand{\id}{\mathrm{Id}}
\newcommand*{\X}{\mathfrak{X}}
\newcommand*{\im}{\mathrm{Im}}
\renewcommand*{\ker}{\operatorname{Ker}}

\newcommand{\abs}[1]{\left\lvert #1 \right\rvert}
\newcommand{\aabs}[1]{\left\lVert #1 \right\rVert}
\newcommand{\ip}[2]{\left\langle #1,#2 \right\rangle}
\DeclareMathOperator{\spt}{spt}

\DeclareMathOperator{\diver}{div}

\newcommand{\hcen}[1]{\hspace*{\fill} #1 \hspace*{\fill}}

\begin{document}


\pagenumbering{gobble}

\topskip0pt
\vspace*{\fill}

\hcen{{\LARGE \bf Keijo M\"onkk\"onen}}
\vspace{25mm}

\begin{center}
{\Huge \bf Integral geometry and unique continuation principles}
\end{center}
\vspace{40mm}

\begin{center}
Esitet\"a\"an Jyv\"askyl\"an yliopiston matemaattis-luonnontieteellisen tiedekunnan 
suostumuksella julkisesti tarkastettavaksi
elokuun 13. p\"aiv\"an\"a 2021 kello 12. 
\end{center}

\vspace{10mm}

\begin{center}
Academic dissertation to be publicly discussed, by permission of 
the Faculty of Mathematics and Science of the University of Jyv\"askyl\"a, 
on August 13, 2021 at 12 o’clock noon. 
\end{center}

\vspace*{\fill}

\hcen{JYV\"ASKYL\"A 2021}
\clearpage


\pagenumbering{roman}

\section*{Foreword}
I wish to thank my supervisor Joonas Ilmavirta for his support and guidance during my years as a PhD student of mathematics. He has helped me to become a little bit more mathematician after my short career in theoretical physics, but not to lose the physicist's way of thinking and calculating without thinking. I also want to thank Mikko Salo for offering me this rare but great opportunity to study inverse problems in his internationally recognized inverse problems group. I want to express my gratitude to the Department of Mathematics and Statistics of University of Jyväskylä for providing me a fruitful working environment in 2017–2021.

I thank Angkana R\"uland for agreeing to be my opponent
at the public examination of my dissertation. I also want to thank the pre-examiners J\"urgen Frikel and Venky Krishnan for their valuable feedback.
I wish to thank my colleagues Jesse Railo and Giovanni Covi for many inspiring moments in mathematics and outside mathematics. Sharing the office room with Jesse has been very beneficial for my research and career, and Giovanni's work on fractional inverse problems convinced me to become interested in that area of mathematics.
I also thank Gunther Uhlmann for collaboration in a very interesting article on higher order fractional Calder\'on problems.
I want to thank all the members in Mikko Salo's inverse problems group, former and present, for the friendly, supporting and encouraging spirit that has always existed in our group.

Finally, I want to thank my family and friends for their support in my long and undirected journey in the academic world.
\\ \bigskip
\bigskip
\bigskip

Jyv\"askyl\"a, \today

Department of Mathematics and Statistics

University of Jyv\"askyl\"a

Keijo M\"onkk\"onen

\clearpage


\section*{List of included articles}
\vspace{5mm}

This dissertation consists of an introductory part and the following seven articles:

\vspace{1em}
\begingroup
\renewcommand{\section}[2]{}

\endgroup










\vspace{3mm}
The author has participated actively in the research of the joint articles~\cite{item:rieszpotential, item:poincareucp, item:covectorpartialdata, item:gunther, item:mixingray, item:polynomialucp}.

\clearpage


\section*{Abstract}
In this thesis we study inverse problems in integral geometry and non-local partial differential equations. We will study these rather different areas of mathematical inverse problems by using the theory of non-local fractional operators. This thesis mainly focuses on proving different kind of unique continuation results of fractional operators which are then used to prove uniqueness results for fractional Calder\'on problems and partial data problems in scalar and vector field tomography.

The introductory part of the thesis contains a general introduction and review of inverse problems arising in medical and seismic imaging. The included articles are divided into three classes which are then presented in their own sections and studied in different levels of detail.

In the articles~\cite{item:rieszpotential, item:poincareucp, item:covectorpartialdata, item:polynomialucp} we consider partial data problems in the X-ray tomography of scalar and vector fields. In the first article~\cite{item:rieszpotential} we prove unique continuation for certain Riesz potentials and apply it to partial data problems of scalar fields. In the second article~\cite{item:poincareucp} we prove unique continuation results for higher order fractional Laplacians which are then used in proving uniqueness for partial data problems of $d$-plane transforms. In the third article~\cite{item:covectorpartialdata} we study partial data problems of vector fields and we prove unique continuation of the normal operator of vector fields which implies uniqueness for the partial data problems.
In the seventh article~\cite{item:polynomialucp} we generalize the unique continuation result of fractional Laplacians proved in~\cite{item:poincareucp} and use it to prove uniqueness for partial data problems of scalar and vector fields, extending the partial data results of the articles~\cite{item:rieszpotential, item:poincareucp, item:covectorpartialdata} to more general cases.

In the articles~\cite{item:poincareucp, item:gunther} we consider higher order fractional Calder\'on problems. In the second article~\cite{item:poincareucp} we use the unique continuation of higher order fractional Laplacians to prove uniqueness for the Calder\'on problem of the higher order fractional (magnetic) Schr\"odinger equation. 
In the fourth article~\cite{item:gunther} we generalize the uniqueness result proved in~\cite{item:poincareucp} to include general lower order local perturbations of the fractional Laplacian.

In the articles~\cite{item:mixingray, item:randers} we consider the travel time tomography problem and its different linearized versions.
In the fifth article~\cite{item:mixingray} we study mixing ray transforms which are generalizations of the geodesic ray transform. We prove solenoidal injectivity results for them in various different cases. 
In the sixth article~\cite{item:randers} we study the boundary rigidity problem on certain non-reversible Finsler manifolds which are also called Randers manifolds. We prove that if the Randers metric consists of a boundary rigid Riemannian metric and a closed 1-form, then the boundary distances determine the Randers metric uniquely up to a natural gauge. 

\clearpage

\begin{otherlanguage}{finnish}
\section*{Tiivistelm\"a}
Tässä väitöskirjassa tutkitaan integraaligeometrian ja epälokaalien osittaisdifferentiaaliyhtälöiden inversio-ongelmia. Näitä melko erilaisia matemaattisia inversio-ongelmia tutkitaan käyttämällä apuna epälokaalien fraktionaalisten operaattoreiden teoriaa. Väitöskirja keskittyy pääosin todistamaan fraktionaalisten operaattoreiden erilaisia yksikäsitteisen jatkon tuloksia, joita käytetään todistaessa yksikäsitteisyyttä fraktionaalisille Calder\'onin ongelmille sekä skalaari- ja vektorikenttien tomografian osittaisen datan ongelmille.

Väitöskirjan johdantokappale sisältää yleisen tason johdatuksen sekä kirjallisuuskatsauksen lääketieteellisessä ja seismisessä kuvantamisessa esiintyviin inversio-ongelmiin. Väitöskirjaan sisällytetyt artikkelit on jaettu kolmeen luokkaan, jotka esitellään omissa kappaleissaan ja joita tarkastellaan yksityiskohtien osalta monella eri tasolla.

Artikkelit~\cite{item:rieszpotential, item:poincareucp, item:covectorpartialdata, item:polynomialucp} käsittelevät osittaisen datan ongelmia skalaari- ja vektorikenttien röntgentomografiassa. Ensimmäisessä artikkelissa~\cite{item:rieszpotential} todistetaan yksikäsitteinen jatko tietyille Rieszin potentiaaleille ja sitä sovelletaan skalaarikenttien osittaisen datan ongelmiin. Toisessa artikkelissa~\cite{item:poincareucp} todistetaan yksikäsitteisen jatkon tuloksia korkeamman kertaluvun fraktionaalisille Laplace-operaattoreille ja niitä käytetään $d$-tasomuunnosten osittaisen datan ongelmien yksikäsitteisyyden todistamisessa.
Kolmannessa artikkelissa~\cite{item:covectorpartialdata} tutkitaan vektorikenttien osittaisen datan ongelmia ja
todistetaan vektorikenttien normaalioperaattorin yksikäsitteinen jatko, josta seuraa yksikäsitteisyys osittaisen datan ongelmille. Seitsemännessä artikkelissa~\cite{item:polynomialucp} yleistetään artikkelissa~\cite{item:poincareucp} todistettu fraktionaalisen Laplace-operaattorin yksikäsitteisen jatkon tulos ja sitä käytetään skalaari- ja vektorikenttien osittaisen datan ongelmien yksikäsitteisyyden todistamisessa, laajentaen artikkeleiden~\cite{item:rieszpotential, item:poincareucp, item:covectorpartialdata} osittaisen datan tuloksia yleisempiin tapauksiin.

Artikkelit~\cite{item:poincareucp, item:gunther} käsittelevät korkeamman kertaluvun fraktionaalisia Calder\'onin ongelmia. Toisessa artikkelissa~\cite{item:poincareucp} käytetään korkeamman kertaluvun fraktionaalisten Laplace-operaattoreiden yksikäsitteistä jatkoa todistaessa yksikäsitteisyyttä korkeamman kertaluvun fraktionaalisen (magneettisen) Schr\"odingerin yhtälön Calder\'onin ongelmalle. Neljännessä artikkelissa~\cite{item:gunther} yleistetään artikkelin~\cite{item:poincareucp} yksikäsitteisyystulos fraktionaalisen Laplace-operaattorin yleisille alempiasteisille lokaaleille perturbaatioille.

Artikkelit~\cite{item:mixingray, item:randers} käsittelevät matka-aikatomografiaa ja sen linearisoituja versioita. Viidennessä artikkelissa~\cite{item:mixingray} tutkitaan sekoitussädemuunnoksia, jotka ovat geodeettisen sädemuunnoksen yleistyksiä. Niille todistetaan solenoidisia injektiivisyystuloksia monissa eri tilanteissa. Kuudennessa artikkelissa~\cite{item:randers} tutkitaan reunajäykkyysongelmaa tietyillä ei-reversiibeleillä Finsler-monistoilla, joita kutsutaan myös Randers-monistoiksi. Artikkelissa todistetaan, että jos Randers-metriikka koostuu reunajäykästä Riemannin metriikasta ja suljetusta 1-muodosta, niin Randers-metriikka määräytyy reunaetäisyyksistään luonnollista mittaa vaille yksikäsitteisesti.

\end{otherlanguage}

\clearpage


\pagenumbering{arabic}

\tableofcontents
\newpage

\section{Introduction}
This thesis is about mathematical inverse problems and the main focus is in proving uniqueness results for different problems arising in tomography. As the title of the thesis suggests, some of the inverse problems appear in integral geometry. However, there are also included inverse problems which do not strictly fit under this category, but they are related to problems in integral geometry via unique continuation principles of non-local operators. 

The inverse problems studied in this thesis can be roughly divided into three classes:
\medskip
\begin{enumerate}[label=(I\arabic*)]
    \item\label{item:traveltimetomographyproblems} The travel time tomography problem and its linearized versions
    \smallskip
    \item\label{item:partialdataproblems}
    Partial data problems in X-ray tomography
    \smallskip
    \item\label{item:fractionalcalderonproblems} Fractional Calder\'on problems.
\end{enumerate}
\medskip
The classes~\ref{item:traveltimetomographyproblems} and~\ref{item:partialdataproblems} belong to integral geometry. In fact, problems in~\ref{item:partialdataproblems} are linearized travel time tomography problems in Euclidean background with partial data. Hence~\ref{item:partialdataproblems} can be seen as a subset of~\ref{item:traveltimetomographyproblems}.
The class~\ref{item:fractionalcalderonproblems} belongs to non-local partial differential equations and at first sight has nothing to do with the classes~\ref{item:traveltimetomographyproblems} and~\ref{item:partialdataproblems}. But there is a way to get from~\ref{item:fractionalcalderonproblems} to~\ref{item:traveltimetomographyproblems}, namely using the ``intermediate step"~\ref{item:partialdataproblems}.

A unifying theme between fractional Calder\'on problems~\ref{item:fractionalcalderonproblems} and partial data problems in X-ray tomography~\ref{item:partialdataproblems} is the use of unique continuation properties of non-local operators in proving uniqueness results. The central operator of this thesis is the fractional Laplace operator~$\fraclaplace$, $s\in (-n/2, \infty)\setminus\Z$, and many of the main theorems of this thesis are unique continuation results of~$\fraclaplace$ or corollaries of them. The unique continuation of~$\fraclaplace$ is used to prove Runge approximation and hence uniqueness for fractional Calder\'on problems. As a special case of fractional Laplacians we have the normal operators of different X-ray transforms whose unique continuation properties are then used to prove uniqueness for various partial data problems arising in the X-ray tomography of scalar and vector fields. 

This introductory part is organized in the following way.
We first discuss in section~\ref{sec:mindmap} how the different articles of this thesis are related to each other. Then we give a gentle introduction to inverse problems and forward problems in section~\ref{sec:inverseanddirect}, and in sections~\ref{sec:xraytomography}--\ref{sec:traveltimetomography} we review the main three classes of inverse problems~\ref{item:traveltimetomographyproblems}--\ref{item:fractionalcalderonproblems} which are studied in this thesis. In sections~\ref{sec:xraypartialdata}--\ref{sec:traveltimefinsler} we go through the main theorems of the included articles. In the beginning of each section we first introduce the inverse problem and give the main results in a general level. We then go through the needed notation in sections~\ref{sec:notationxray}--\ref{sec:notationtraveltime} before giving the main theorems with all technical details in sections~\ref{sec:mainresultsxraypartial}--\ref{sec:mainresultstraveltime}. Section~\ref{sec:traveltimefinsler} can be read independently of sections~\ref{sec:xraypartialdata} and~\ref{sec:fractionalcalderon}. Section~\ref{sec:fractionalcalderon} can also be read independently of section~\ref{sec:xraypartialdata} if one first goes through the notation in section~\ref{sec:notationxray}.

\subsection{On the articles of this thesis}
\label{sec:mindmap}
\begin{figure}[htp]
\centering
\captionsetup{margin=0cm}
\includegraphics[height=6.85cm]{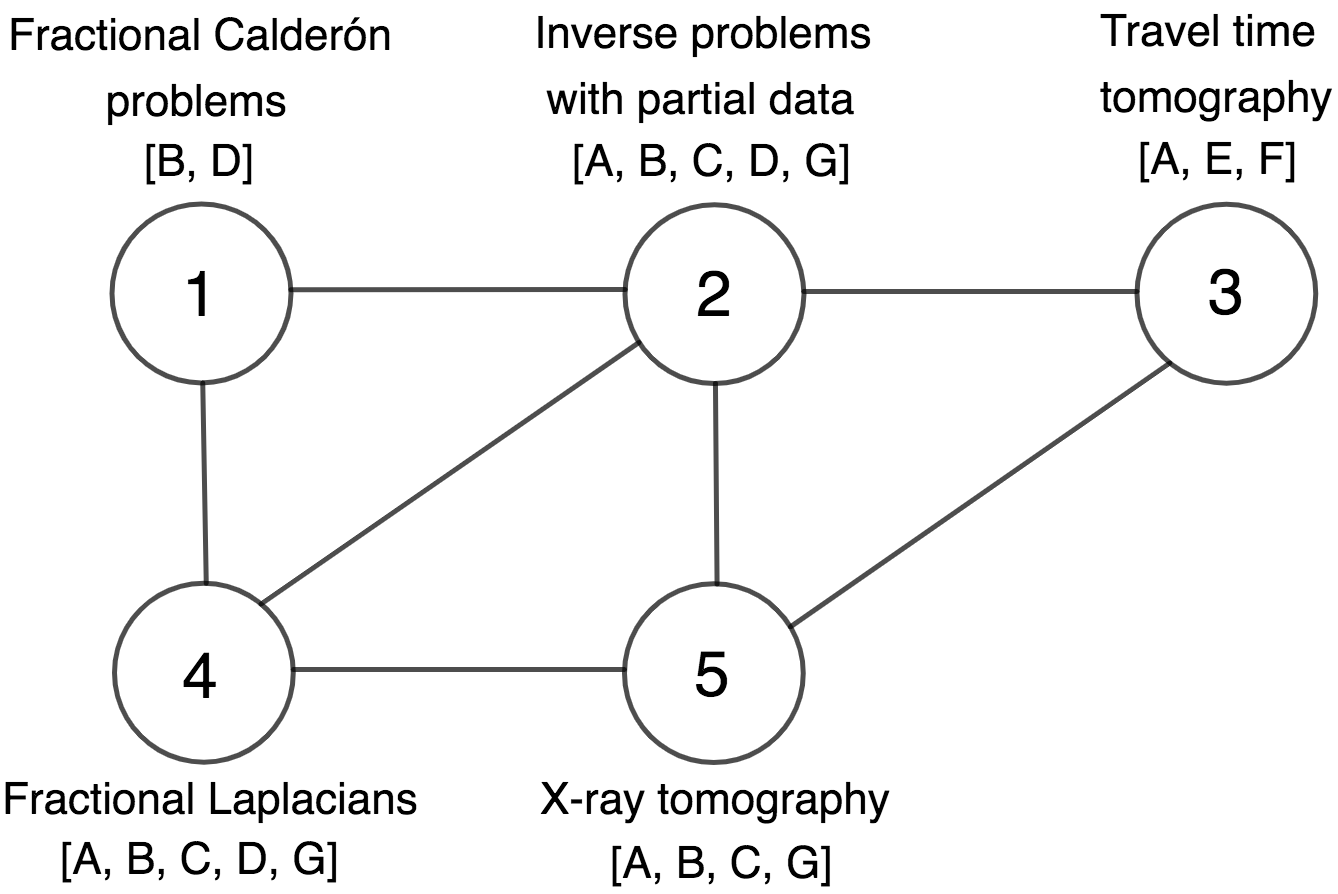}
\caption{A graph illustrating the relation between the different articles of this thesis.}
\label{fig:mindmap}
\end{figure}

In figure~\ref{fig:mindmap} we have illustrated the connection between the different articles of this thesis. In most of the articles we study inverse problems with partial data: these include fractional Calder\'on problems (articles~\cite{item:poincareucp, item:gunther}), X-ray tomography with partial data (articles~\cite{item:rieszpotential, item:poincareucp, item:covectorpartialdata, item:polynomialucp}) and linearized travel time tomography with ``half-local" data (article~\cite{item:rieszpotential}). Unique continuation of fractional Laplacians has a crucial role in proving uniqueness for partial data problems studied in this thesis. Fractional Laplacians arise in fractional Calder\'on problems and also in X-ray tomography in the form of different normal operators. Problems in X-ray tomography in turn can be seen as linearized travel time tomography problems in Euclidean background.

In fractional Calder\'on problems the task is to recover the potential (and more generally a perturbation) of the fractional Schr\"odinger equation in a bounded domain by doing measurements in the exterior of the domain. These problems are studied in the articles~\cite{item:poincareucp, item:gunther} and treated in section~\ref{sec:fractionalcalderon}. In X-ray tomography we want to determine a scalar field (or a vector field) when we know its integrals over lines which intersect a given nonempty open set. This is studied in the articles~\cite{item:rieszpotential, item:poincareucp, item:covectorpartialdata, item:polynomialucp} and treated in section~\ref{sec:xraypartialdata}. In travel time tomography one wants to recover the speed of sound (and more generally Riemannian metric or Finsler norm) by measuring travel times (geodesic distances) on the boundary of a compact manifold. This problem and its linearized versions (the geodesic ray transform and its generalizations) are studied in the articles~\cite{item:rieszpotential, item:mixingray, item:randers} and treated in sections~\ref{sec:xraypartialdata} and~\ref{sec:traveltimefinsler}. 

A remark from the point of view of graph theory: the (connected) graph presented in figure~\ref{fig:mindmap} has a Hamilton cycle, i.e. a closed walk such that every vertex is visited exactly once. The graph also has an Euler trail, i.e. a walk such that every edge is traversed exactly once. Formally, this ``proves" that the articles of this thesis are closely related to each other. However, the graph does not admit an Euler tour (a closed walk which is an Euler trail) since not every vertex has even degree~\cite{BM-graph-theory-applications}.

\subsection{Inverse problems and forward problems}
\label{sec:inverseanddirect}
Inverse problems are practical or abstract problems which arise for example in medical and seismic imaging~\cite{HR-inverse-problems, ISA-inverse-source-problems, KI-introduction-to-inverse-problems, MS-linear-nonlinear-inverse-problems, NA-mathematics-computerized-tomography, NA-mathematical-methods-image-reconstruction, NN-intro-inverse-problems, RK-radon-transform-an-local-tomography, UH-inverse-problems-seeing-the-unseen}. Inverse problems are often encountered when making indirect measurements. In such situations we have an object we cannot or do not want to access by invasive methods. In medical imaging the object can be a patient we want to study without doing surgical operations, and in seismic imaging the object can be the planet Earth whose deep interior we cannot reach by any practical means. The common task in both cases is that one wants to deduce the interior features of some object by making measurements on the boundary or in the exterior of the object. Usually we have some physical model which tells us how the interior properties of the object affect the measurements we make on the boundary or in the exterior. The goal is to use this physical model to deduce the interior properties of the object from the boundary or exterior measurements. The boundary and exterior measurements are often called just data.

Inverse problems are opposite to what we call direct problems or forward problems. Let us consider an example from X-ray tomography to illustrate the difference. In X-ray tomography one shoots X-rays through an object and studies the attenuation pattern of the X-rays. The attenuation of the X-rays is determined by the interior properties (the position-dependent attenuation coefficient) of the object. In the direct problem one knows the attenuation of the object and wants to determine the attenuation pattern of the X-rays. When the initial intensity of the X-rays is known, then one can easily calculate the final intensity of the X-rays by using a simple physical model~\cite{NA-mathematics-computerized-tomography}. Roughly saying, the direct problem corresponds to putting values for parameters in an equation and computing the result.

Inverse problems are much harder since they ``operate" in the opposite direction. For example, in medical imaging one wants to determine the attenuation of the object instead of the attenuation pattern of the X-rays which can be easily measured. Since one also can control the initial intensity of the X-rays we have indirect information about the attenuation, i.e. we know the total attenuation of the X-rays and want to determine the attenuation of the object from that data. It turns out that the total attenuation corresponds to the integrals of the attenuation function along lines which intersect the object~\cite{NA-mathematics-computerized-tomography}. The inverse problem is to invert this integral transform which is also called the X-ray transform. The inversion of the X-ray transform is a much harder task than solving the forward problem where we already know the interior features of the object (the attenuation) and just have to calculate the end result (the final intensity of the X-rays).

Uniqueness, stability and reconstruction are important properties in the study of inverse problems. Uniqueness means that the inverse problem has a unique solution. In other words, if two objects produce the same boundary or exterior data, then they must have the same interior features. Reconstruction means that there is some way (e.g. an algorithm or formula) so that one can compute the desired physical quantity related to the interior properties of the object from the boundary or exterior data. Stability is related to how much measurement errors affect uniqueness or reconstruction. Since in practice there is always some noise in measurements, stability is important in showing that the reconstructed quantity is not too far away from the true value of that quantity. These three properties are not independent of each other since uniqueness usually follows from reconstruction and stability.

Uniqueness and stability have a connection to Hadamard's formulation of a well-posed problem~\cite{HA-hadamard-book, HA-hadamard-book2}. A mathematical problem related to a physical phenomenon is called well-posed, if the problem has unique solution which is stable with respect to the measured data (the solution depends continuously on the data)~\cite{HR-inverse-problems, KI-introduction-to-inverse-problems, MS-linear-nonlinear-inverse-problems, NA-mathematics-computerized-tomography, RK-radon-transform-an-local-tomography}. 
If the solution fails to exist, the solution is not unique or the solution does not depend continuously on the data, the problem is said to be ill-posed.
Forward problems are often well-posed, but inverse problems tend to be ill-posed. Usually the reason for ill-posedness of inverse problems is that they lack stability which causes difficulties in numerical reconstruction~\cite{KI-introduction-to-inverse-problems, KRS-instability-of-inverse-problems, MS-linear-nonlinear-inverse-problems, NN-intro-inverse-problems, RK-radon-transform-an-local-tomography}. 

In this thesis we mainly focus on uniqueness, i.e. in most of our theorems we show that the inverse problem has unique solution.
Even if we do not get stability or a reconstruction formula for the problem, uniqueness is important in practical applications. Uniqueness for example increases the reliability of the results obtained in X-ray tomography when we only have a ``small amount" of measurement data available.

\subsection{X-ray tomography of scalar and vector fields}
\label{sec:xraytomography}
\subsubsection{X-ray tomography of scalar fields}
X-ray tomography is a commonly used method in medical imaging to study interiors of objects. The main goal is to determine the attenuation of the object when one knows the initial and final intensities of X-rays, i.e. the total attenuation of the X-rays. The attenuation can be modelled as a scalar function~$f\colon\R^n\to\R$. If the rays propagate parallel to $x$-axis, then the intensity~$I$ of the X-rays satisfies the differential equation $I'(x)=-f(x)I(x)$ and the total attenuation corresponds to the line integral~\cite{NA-mathematics-computerized-tomography}
\begin{equation}
\label{eq:totalattenuation}
\ln\bigg(\frac{I_0}{I_1}\bigg)=\int_\gamma f\der s
\end{equation}
where~$I_0$ is the initial intensity and~$I_1$ is the final intensity of the X-rays, and~$\gamma$ is a line along which the X-ray beam propagates. The inverse problem in X-ray tomography is to solve~$f$ in equation~\eqref{eq:totalattenuation} using different lines~$\gamma$ when the left-hand side of the equation is known.

The previous discussion motivates us to define the operator~$\xrt_0$ as
\begin{equation}
\label{eq:xraytransformscalar}
\xrt_0 f(\gamma)=\int_\gamma f\der s
\end{equation}
where~$\gamma$ is a line in~$\R^n$ and~$f\colon\R^n\rightarrow \R$ is a scalar field. The operator~$\xrt_0$ is called the X-ray transform of scalar fields and in two dimensions~$\xrt_0$ is also known as the Radon transform. 
The inverse problem is to invert the operator~$\xrt_0$ in equation~\eqref{eq:xraytransformscalar} and it was first studied by Johann Radon~\cite{RA-radon-transform}. Theoretical and practical applications to computerized tomography were studied by Cormack and Hounsfield~\cite{CO-ct-imaging, CO-ct-imaging2, NA-mathematics-computerized-tomography}. There are formulas for the inversion of~$\xrt_0$ some of which involve the normal operator~$\no_0$ of the X-ray transform~\cite{HE:integral-geometry-radon-transforms, NA-mathematics-computerized-tomography, RK-radon-transform-an-local-tomography, SU:microlocal-analysis-integral-geometry}. The normal operator~$\no_0=\xrt_0^*\xrt_0$ is defined as first applying the X-ray transform and then back-projecting~$\xrt_0 f$ from the space of all lines to a function in~$\R^n$ using the adjoint~$\xrt_0^*$ of the X-ray transform. Hence~$\no_0$ is a useful auxiliary operator which maps functions on~$\R^n$ to functions on~$\R^n$ and one can study the X-ray transform~$\xrt_0$ using its normal operator~$\no_0$.

The inversion formulas for~$\xrt_0$ assume that we know the integrals of~$f$ over all lines in~$\R^n$. In practical applications we only have access to a small subset of lines, and in that case we have a partial data problem. One such partial data problem is to uniquely determine~$f$ everywhere in~$\R^n$ from its X-ray data on all lines intersecting a given open set~$V\subset\R^n$. The integrals alone cannot determine~$f$ uniquely and one has to make additional assumptions~\cite{KEQ-wavelet-methods-ROI-tomography, NA-mathematics-computerized-tomography}. The partial data problem has unique solution, if $f|_V$ vanishes~\cite{CNDK-solving-interior-problem-ct-with-apriori-knowledge, KEQ-wavelet-methods-ROI-tomography}, $f|_V$ is piecewise constant or piecewise polynomial~\cite{KEQ-wavelet-methods-ROI-tomography, YYJW-high-order-TV-minimization} or if $f|_V$ is real analytic~\cite{KKW-stability-of-interior-problems}. A complementary partial data result is the Helgason support theorem where one has access to lines which do not intersect a given compact and convex set and the problem is to determine the scalar field uniquely outside that set~\cite{HE:integral-geometry-radon-transforms}. Partial data problems are in general much harder to treat than problems with full data because the reconstruction is not stable anymore and there can be artefacts in the images even if the problem admits a unique solution. In such cases we have ``invisible singularities"~\cite{KQ-microlocal-analysis-in-tomography, KLM-on-local-tomography, NA-mathematics-computerized-tomography, QU-singularities-x-ray-transform-limited-data, QU-artifacts-and-singularities-limited-tomography}.

We can generalize the transform~$\xrt_0$ from lines to affine $d$-dimensional planes where $0<d<n$. The $d$-plane transform~$\dplane$ is defined as~\cite{HE:integral-geometry-radon-transforms}
\begin{equation}
\label{eq:dplanetransform}
\dplane f(A)=\int_{A}f(x)\der m(x)
\end{equation}
where~$A$ is an affine $d$-dimensional plane, ~$m$ is the $d$-dimensional Hausdorff measure and~$f\colon\R^n\rightarrow\R$ is a scalar field. The case $d=1$ corresponds to the X-ray transform~$\xrt_0$ and the case $d=n-1$ is often called the Radon transform which coincides with the X-ray transform in two dimensions~\cite{HE:integral-geometry-radon-transforms, NA-mathematics-computerized-tomography, RK-radon-transform-an-local-tomography}. As before, the inverse problem is to invert the transform~$\dplane$ in equation~\eqref{eq:dplanetransform}. There is an inversion formula in terms of the normal operator~$\nod$ of the $d$-plane transform which is defined in a similar way as in the case of the X-ray transform~\cite{HE:integral-geometry-radon-transforms}. One example of partial data results for $d$-plane transforms is the Helgason support theorem where one knows the integrals of the scalar field over all $d$-planes which do not intersect a given compact and convex set~\cite{HE:integral-geometry-radon-transforms}.

The $d$-plane transform (also called the $k$-plane transform in some works) has been extensively studied after the pioneering work by Fuglede~\cite{FU-d-plane} and Helgason~\cite{HE-d-plane-transforms}. See for example~\cite{A11, GO-range-of-d-plane-transform, HE:integral-geometry-radon-transforms, I15, KU-d-plane, RA-periodic-radon-transform} and the works by Rubin~\cite{RU-k-plane-transforms4, RU-k-plane-transforms5, RU-k-plane-transforms3, RU-k-plane-transform, RU-k-plane-transforms2}.

\subsubsection{X-ray tomography of vector fields}
X-ray tomography is also used in the imaging of moving fluids which is based on Doppler backscattering or acoustic travel time measurements. 
If~$h\colon\R^n\rightarrow\R^n$ is a vector field which represents the flow field of a moving fluid, then after a linearization procedure one ends up studying the transform~\cite{NO-tomographic-recostruction-of-vector-fields, NO-unique-tomographic-reconstruction-doppler}
\begin{equation}
\label{eq:xraytransformvector}
\xrt_1 h(\gamma)=\int_\gamma h\cdot\der\overline{s}.
\end{equation}
The operator~$\xrt_1$ is called the X-ray transform of vector fields and it has applications for example in medical ultrasound imaging~\cite{JASESPL-ultrasound-doppler-tomography, JUH-principles-of-doppler-tomography, SCHU-importance-of-vector-field-tomography, SSLP-doppler-tomography-vector-fields}.
The inverse problem is to invert the operator~$\xrt_1$ in equation~\eqref{eq:xraytransformvector}.

Unlike in the scalar case we have a natural gauge: the gradients of scalar fields which vanish at infinity are always in the kernel of~$\xrt_1$. For this reason one can determine the vector field~$h$ only up to potential fields from its X-ray transform, i.e. one can only determine the solenoidal part~$\sol{h}$ in the Helmholtz decomposition $h=\sol{h}+\nabla\phi$ where $\sol{h}\colon\R^n\rightarrow\R^n$ is a vector field such that $\diver \sol{h}=0$ and~$\phi$ is a scalar field~\cite{SCHU-importance-of-vector-field-tomography, SHA-integral-geometry-tensor-fields, SU:microlocal-analysis-integral-geometry}. The solenoidal part can be uniquely determined from the full X-ray data and there is an inversion formula in terms of the normal operator~$\no_1=\xrt_1^*\xrt_1$ of the X-ray transform of vector fields where $\xrt_1^*$ is the adjoint operator (or back-projection)~\cite{JUH-principles-of-doppler-tomography, NO-tomographic-recostruction-of-vector-fields, SHA-integral-geometry-tensor-fields, SSLP-doppler-tomography-vector-fields, SU:microlocal-analysis-integral-geometry}.

Like in the scalar case, one can also study X-ray tomography of vector fields with partial data. The main goal in such problems is to determine the solenoidal part of the vector field from its partial X-ray data. Examples of such partial data results include cases where one knows the integrals of the vector field over lines which intersect a certain type of curve~\cite{DEN-inversion-of-3d-tensor-fields, RA-microlocal-analysis-doppler-transform, VER-integral-geometry-symmetric-tensor-incomplete} or which are parallel to a finite set of planes~\cite{JUH-principles-of-doppler-tomography, SCHU-3d-doppler-transform-reconstruction-and-kernels, SHA-vector-tomography-incomplete-data}. There is also a vectorial version of the Helgason support theorem where one knows the integrals of the vector field over all lines not intersecting a given convex and compact set~\cite{SU:microlocal-analysis-integral-geometry}.
\\ \\
\subsection{Electrical impedance tomography and its non-local versions}
\label{sec:eit}
\subsubsection{The Calder\'on problem}
Electrical impedance tomography (EIT) is an imaging method which has applications in geophysics and medical imaging~\cite{KI-introduction-to-inverse-problems, MS-linear-nonlinear-inverse-problems, UH-inverse-problems-seeing-the-unseen}. EIT is based on the conductivity equation and the inverse problem is known as the Calder\'on problem. In the Calder\'on problem we have an object whose electrical properties we want to deduce by making boundary measurements. In particular, we want to determine the conductivity inside the object by applying voltages on the boundary and measuring the induced currents on the boundary which depend on the electrical properties of the interior of the object.

We can model the object as a bounded domain $\Omega\subset\R^n$ with sufficiently regular boundary~$\partial\Omega$. The conductivity equation is~\cite{UH-inverse-problems-seeing-the-unseen}
\begin{align}\label{eq:conductivityequation} \left\{\begin{array}{rl}
        \nabla\cdot(\eta\nabla u)&=0 \;\;\text{in}\;  \Omega\\
        u|_{\partial\Omega} &=f
        \end{array}\right.\end{align}
where~$f$ is the potential on the boundary, $u$ is the induced potential in~$\Omega$ and~$\eta$ is the electrical conductivity of~$\Omega$ which is assumed to be sufficiently smooth positive function. The measurements are encoded in the Dirichlet-to-Neumann (DN) map~$\Lambda_\eta$ which tells how the electrical properties of the interior induce normal currents on the boundary when one applies the voltage~$f$ on the boundary. More specifically, one can write $\Lambda_\eta f=(\eta\partial_\nu u)|_{\partial\Omega}$ where~$\nu$ is the outer unit normal on $\partial\Omega$. The inverse problem is to determine the conductivity~$\eta$ in equation~\eqref{eq:conductivityequation} by applying different boundary values~$f$ (voltages) and measuring the induced currents~$\Lambda_\eta f$. In particular, the uniqueness problem is the following: if $\Lambda_{\eta_1}f=\Lambda_{\eta_2}f$ for all boundary values~$f$, does it follow that $\eta_1=\eta_2$? This problem was first studied mathematically by Alberto Calder\'on and the inverse problem is therefore known as the Calder\'on problem~\cite{CA-calderon-problem}.

Using the substitution $\tilde{u}=\sqrt{\eta}u$ one can convert the conductivity equation~\eqref{eq:conductivityequation} to the following Schr\"odinger equation~\cite{NA-calderon-problem, SU-cgo-conductivity, UH-inverse-problems-seeing-the-unseen}
\begin{align}\label{eq:scrodingercalderon} \left\{\begin{array}{rl}
        (-\Delta+q)\tilde{u}&=0  \;\;\text{in}\;  \Omega\\
        \tilde{u}|_{\partial\Omega} &=\tilde{f}.
        \end{array}\right.\end{align}
Here~$q=(\Delta\sqrt{\eta})/\sqrt{\eta}$ now corresponds to the electric potential in~$\Omega$ and $\tilde{f}=\sqrt{\eta}f$. The DN map~$\Lambda_q$ for equation~\eqref{eq:scrodingercalderon} can be written as $\Lambda_q \tilde{f}=\partial_{\nu}\tilde{u}|_{\partial\Omega}$ assuming~$\partial\Omega$ is regular enough. The interpretation of the DN map is as in the conductivity equation: the DN map tells how the applied voltage on the boundary induces normal currents on the boundary via the electrical properties of the interior of the object. 
The inverse problem now is to determine the potential~$q$ in equation~\eqref{eq:scrodingercalderon} by applying different boundary values~$\tilde{f}$ (voltages) and measuring the induced currents~$\Lambda_q \tilde{f}$. The uniqueness problem is as for the conductivity equation: if $\Lambda_{q_1}\tilde{f}=\Lambda_{q_2}\tilde{f}$ for all boundary values~$\tilde{f}$, does it follow that $q_1=q_2$? One standard tool in proving uniqueness for the Calder\'on problem of the conductivity equation~\eqref{eq:conductivityequation} and Schr\"odinger equation~\eqref{eq:scrodingercalderon} is the construction of complex geometrical optics solutions~\cite{AP-calderon-two-dimensions, CA-calderon-problem, SU-cgo-conductivity2, SU-cgo-conductivity, UH-inverse-problems-seeing-the-unseen}.  

\subsubsection{The fractional Calder\'on problem}
One can study the non-local version of the Schr\"odinger equation~\eqref{eq:scrodingercalderon} as follows. One replaces the Laplacian~$-\Delta$ with the fractional Laplacian~$\fraclaplace$ which is the pseudodifferential operator
\begin{equation}
\fraclaplace u=\ifourier(\abs{\cdot}^{2s}\hat{u}), \quad s\in (-n/2, \infty)\setminus\Z.
\end{equation}
The fractional Laplacian is a non-local operator in contrast to the ordinary Laplacian: the value~$\fraclaplace u(x)$ depends on the values of~$u$ everywhere in~$\R^n$ while $-\Delta u(x)$ depends only on the values of~$u$ in a small neighborhood of $x\in\R^n$. For example, the normal operator of the X-ray transform~$\no_0$ (and more generally the normal operator of the $d$-plane transform~$\nod$) is the fractional Laplacian $(-\Delta)^{-1/2}$ (more generally $(-\Delta)^{-d/2}$) up to a constant factor. In addition to integral geometry fractional Laplacians arise also in non-local diffusion~\cite{BV-nonlocal-diffusion-applications, DGLZ2012, GSU-calderon-problem-fractional-schrodinger} and in fractional quantum mechanics~\cite{LA-fractional-quantum-mechanics, LA-fractional-quantum-mechanics-book}.

Replacing~$-\Delta$ with~$\fraclaplace$ where $s\in (0, 1)$ we obtain the fractional Schr\"odinger equation introduced in~\cite{GSU-calderon-problem-fractional-schrodinger}
\begin{align}\label{eq:fractionalschrodingerequation} \left\{\begin{array}{rl}
        (\fraclaplace+q)u&=0  \;\;\text{in}\;  \Omega\\
        u|_{\Omega_e} &=f
        \end{array}\right.\end{align}
where~$\Omega_e=\R^n\setminus\overline{\Omega}$ is the exterior of the bounded domain~$\Omega\subset\R^n$. For such non-local equation~\eqref{eq:fractionalschrodingerequation} it is more natural to consider exterior values $u|_{\Omega_e}=f$ instead of boundary values. The DN map~$\Lambda_q$ maps the ``non-local voltage"~$f$ to a non-local version of the normal current~\cite{GSU-calderon-problem-fractional-schrodinger}: under stronger assumptions one can write $\Lambda_q f=\fraclaplace u|_{\Omega_e}$. In the fractional Calder\'on problem one wants to determine the potential~$q$ in equation~\eqref{eq:fractionalschrodingerequation} by applying different exterior values~$f$ and measuring the induced ``exterior currents"~$\Lambda_q f$. The uniqueness problem is similar as in the local case: if $\Lambda_{q_1}f=\Lambda_{q_2}f$ for all exterior values~$f$, does it follow that $q_1=q_2$? The fractional Calder\'on problem for equation~\eqref{eq:fractionalschrodingerequation} was first studied by Ghosh, Salo and Uhlmann~\cite{GSU-calderon-problem-fractional-schrodinger}.

In fractional Calder\'on problems instead of constructing complex geometrical optics solutions one can exploit the non-locality of the equation and especially the non-local behaviour of the operator~$\fraclaplace$. One has the following unique continuation property of fractional Laplacians~\cite{GSU-calderon-problem-fractional-schrodinger}: if $s\in (0, 1)$ and $\fraclaplace u|_V=u|_V=0$ for some nonempty open set $V\subset\R^n$, then $u=0$. Clearly such property cannot hold for local operators such as~$-\Delta$. The unique continuation of~$\fraclaplace$ is in essential role in proving uniqueness for fractional Calder\'on problems~\cite{BGU-lower-order-nonlocal-perturbations, CLR18, CO-magnetic-fractional-schrodinger, GSU-calderon-problem-fractional-schrodinger}. 

After the seminal work~\cite{GSU-calderon-problem-fractional-schrodinger} there have been numerous results for different variants of the fractional Calder\'on problem: these include stability and instability results~\cite{RS18, RS-fractional-calderon-low-regularity-stability}, uniqueness under single measurement~\cite{GRSU-fractional-calderon-single-measurement}, magnetic versions of the fractional Schr\"odinger equation~\cite{CO-magnetic-fractional-schrodinger, LILI-semilinear-magnetic, LI-fractional-magnetic, LILI-fractional-magnetic-calderon}, lower order local and non-local perturbations~\cite{BGU-lower-order-nonlocal-perturbations, CLR18}, semilinear equations~\cite{LL19, LL-fractional-semilinear-problems}, fractional conductivity and heat equations~\cite{Co18, LLR19,ruland2019quantitative} and equations arising from a non-local Schr\"odinger-type elliptic operator~\cite{CLL19, GLX-calderon-nonlocal-elliptic-operators}.

\subsection{Travel time tomography and its linearization}
\label{sec:traveltimetomography}
\subsubsection{The boundary rigidity problem}
In seismic travel time tomography the objective is to study the interior properties of the Earth by measuring travel times of seismic waves on the surface of the Earth~\cite{CE-seismic-ray-theory, HR-inverse-problems, SHE-introduction-to-seismology, SUVZ-travel-time-tomography}. It is impossible to access the deep interior of the Earth by any practical means and the only way to obtain information is by doing indirect measurements on the surface. The travel times of seismic waves depend on the speed of sound in the medium where the wave propagates. Therefore the travel times contain indirect information about the physical properties of the Earth.

The Earth can be modelled as a three-dimensional compact manifold~$M$ with boundary~$\partial M$ (e.g. a closed ball). Assuming that the medium is isotropic the speed of sound depends only on position and it becomes a positive scalar function $c\colon M\to (0, \infty)$. The travel time of a seismic wave or ray can be expressed as the line integral~\cite{CE-seismic-ray-theory}
\begin{equation}
\label{eq:traveltime}
T=\int_\gamma\frac{\der s}{c}
\end{equation}
where~$\gamma$ is the ray path. The travel time tomography problem or inverse kinematic problem is to solve the speed of sound~$c$ in equation~\eqref{eq:traveltime} when the travel times~$T$ measured on the surface are known.

The travel time tomography problem was studied first in 1900s by Herglotz, Wiechert and Zoeppritz~\cite{HE-inverse-kinematic-problem, WZ-kinematic-problem}. They solved the problem assuming that the speed of sound is radial $c=c(r)$ and satisfies the Herglotz condition
\begin{equation}
\label{eq:herglotzcondition}
\frac{\der}{\der r}\bigg(\frac{r}{c(r)}\bigg)>0.
\end{equation}
Under these assumptions the solution reduces to the inversion of an Abel-type integral transform~\cite{NOW-herglotz-abel-transform, SHE-introduction-to-seismology}.
The Herglotz condition~\eqref{eq:herglotzcondition} is equivalent to the condition that the travel times in equation~\eqref{eq:traveltime} are finite~\cite{deI:abel-transforms-x-ray-tomography}. In geometrical terms, the Herglotz condition~\eqref{eq:herglotzcondition} means that one can foliate the manifold~$M$ with strictly convex hypersurfaces (i.e. spheres)~\cite{SUVZ-travel-time-tomography}.

The travel time tomography problem can be formulated in a more geometrical way. The speed of sound~$c$ determines the Riemannian metric $g_c=c^{-2}(x)e$ where~$e$ is the Euclidean metric. By Fermat's principle the rays propagate along geodesics of the metric~$g_c$ and the travel times correspond to lengths of these geodesics~\cite{CE-seismic-ray-theory}. 
The inverse problem is to determine the scalar function~$c$, or equivalently the metric~$g_c$, from the lengths of all geodesics connecting points on the boundary~$\partial M$. One sees that the problem is highly non-linear since the geodesics depend on the function~$c$ (or the metric~$g_c$).

One can study the above geometric problem in a more general case: if~$g$ is a Riemannian metric, determine~$g$ from the distances between boundary points (boundary distances) given by~$g$. This geometric inverse problem is known as the boundary rigidity problem~\cite{SUVZ-travel-time-tomography}. In particular, one problem of interest is the uniqueness problem: if two Riemannian metrics~$g_1$ and~$g_2$ give the same boundary distances, does it follow that $g_1=g_2$? The answer is no in general since there is a gauge: if $g_2=\Psi^*g_1$ where $\Psi\colon M\rightarrow M$ is a diffeomorphism which is identity on the boundary, then~$g_1$ and~$g_2$ give the same boundary distances~\cite{SUVZ-travel-time-tomography}. Hence without further restrictions one can determine the metric only up to a boundary preserving diffeomorphism.

The boundary rigidity problem is a difficult non-linear inverse problem and it has been solved only in certain special cases where the manifold admits strictly convex foliation~\cite{SUV-boundary-rigidity-foliation, SUVZ-travel-time-tomography} or the manifold is known to be simple (a generalization of a Euclidean ball). Boundary rigidity holds for simple subspaces of Euclidean space~\cite{GRO-filling-riemannian-manifolds} and simple subspaces of symmetric spaces of constant negative curvature~\cite{BCG-rigidity}.  In two dimensions examples include simple subspaces of the open hemisphere~\cite{MI-open-hemisphere} and simple spaces of negative curvature~\cite{CRO-rigidity-negative-curvature}.
If the Riemannian metrics on a compact simple Riemannian manifold are in the same conformal class, then the distances between boundary points determine the metric uniquely, i.e. the diffeomorphism~$\Psi$ becomes identity in this case~\cite{CRO-rigidity-conformal, MU-reconstruction-problem-two-dimensional, SUVZ-travel-time-tomography}. In general, compact simple Riemannian manifolds are known to be boundary rigid in two dimensions~\cite{PU-simple-manifolds-boundary-rigidity}, but it is conjectured that boundary rigidity holds for compact simple Riemannian manifolds of any dimension~\cite{RE-boundary-rigidity-conjecture}.

In the travel time tomography problem one usually assumes that the speed of sound~$c$ is isotropic, i.e. it only depends on position. However, anisotropies have been observed in the shallow crust, upper mantle and inner core of the Earth~\cite{CRE-anisotropy-inner-core, DA-PREM-model, SHE-introduction-to-seismology}. Therefore it is reasonable to consider~$c$ as a function on the tangent bundle $c\colon TM\rightarrow (0, \infty)$ so that the dependence on the direction of propagation can be taken into account. If the sound speed is anisotropic, then the seismic rays propagate along geodesics of a Finsler norm and we need Finsler geometry to treat the anisotropies~\cite{ABS-seismic-rays-as-finsler-geodesics, YN-finsler-seismic-ray-path}. The travel time tomography problem can then be expressed as a boundary rigidity problem on Finsler manifolds where the fiberwise inner product depends not only on position but also on direction.

The boundary rigidity problem is much harder in the Finslerian case since there are non-isometric Finsler norms which give the same boundary distances~\cite{BI-boundary-rigidity, CNV-finsler-deformations, CD-length-spectrum-orbifold, IVA-finsler-monotonicity}. This means that in general Finsler norms are not rigid in the same way as Riemannian metrics.
However, some rigidity results are known in certain special cases. Projectively flat Finsler norms on compact convex domains of~$\R^2$ are uniquely determined by their boundary distances~\cite{ALE-convex-rigidity-plane, AM-pseudo-metrics-on-the-plane, KO-boundary-rigidity-projective-metrics}. When we restrict ourselves to Finsler norms which are relevant in seismology, we can expect more rigidity: one can use the collection of boundary distance maps to determine the differential and topological structures of Finsler manifolds~\cite{deILS-finsler-boundary-distance-map}, and the broken scattering relation determines the isometry class of reversible Finsler manifolds which admit strictly convex foliation~\cite{deILS-broken-scattering-rigidity}.

\subsubsection{Linearized versions of the boundary rigidity problem}
Let us study the linearization of the boundary rigidity problem. Let $\epsilon>0$ and $s\in (-\epsilon, \epsilon)$. Assume that~$g^s$ is a family of Riemannian metrics which all give the same boundary distances where~$g^0$ corresponds to a known ``background metric". When we linearize the boundary rigidity problem, we calculate the derivative of the boundary distances at $s=0$. Since these distances do not depend on the parameter~$s$ we obtain~\cite{SHA-integral-geometry-tensor-fields}
\begin{equation}
0=\int_a^b \frac{\partial g^s_{ij}(\gamma_0(t))}{\partial s}\bigg|_{s=0}\dot{\gamma}_0^i(t)\dot{\gamma}_0^j(t)\der t
\end{equation}
where $\gamma_0\colon [a, b]\to M$ is a geodesic of the base manifold $(M, g^0)$ connecting two boundary points. If the variations~$g^s$ are conformal, i.e. $g^s=f_s g_0$ where~$f_s\colon M\rightarrow\R$ is a family of positive scalar functions such that $f_0=1$, then the linearization leads to

\begin{equation}
0=\int_a^b \frac{\partial f_s(\gamma_0(t))}{\partial s}\bigg|_{s=0}\der t.
\end{equation}

The previous observations motivate us to study the kernel of the geodesic ray transform of symmetric $m$-tensor fields where $m\geq 0$.
The geodesic ray transform of a scalar field~$f\colon M\rightarrow\R$ (or 0-tensor field) on a Riemannian manifold $(M, g)$ is defined as
\begin{equation}
\label{eq:geodesicraytransformscalar}
\geod_0 f(\gamma)=\int_{\tau_\gamma^-}^{\tau_\gamma^+}f(\gamma (t))\der t
\end{equation}
where $\gamma\colon [\tau_\gamma^-, \tau_\gamma^+]\to M$ is a geodesic defined on the maximal interval $[\tau_\gamma^-, \tau_\gamma^+]$ which can be finite or infinite.
More generally, the geodesic ray transform of a symmetric (covariant) $m$-tensor field~$h$ is ($m\geq 1$)
\begin{equation}
\label{eq:geodesicraytransform}
\geod_m h(\gamma)=\int_{\tau_\gamma^-}^{\tau_\gamma^+}h_{i_1\dotso i_m}(\gamma (t))\dot{\gamma}^{i_1}(t)\cdots\dot{\gamma}^{i_m}(t) \der t
\end{equation}
where $h_{i_1\dotso i_m}(x)$ are the components of the $m$-tensor field~$h$ in local coordinates and we have used the Einstein summation convention (repeated indices which appear both as a subscript and superscript are implicitly summed over).
The geodesic ray transform~$\geod_m$ can be seen as a generalization of the Euclidean X-ray transform since in Euclidean space geodesics are lines. However, Funk studied the geodesic ray transform of scalar fields on the sphere $S^2\subset\R^3$ (also known as the Funk transform) before Radon introduced the Euclidean X-ray transform or Radon transform~\cite{FU-geodesic-ray-transform, FU-geodesic-ray-transform2, HE:integral-geometry-radon-transforms}.

The inverse problem in geodesic ray tomography is to determine the $m$-tensor field~$h$ (or the scalar field~$f$) from its integrals along geodesics, i.e. we want to invert the operator~$\geod_m$ in equation~\eqref{eq:geodesicraytransform} (or in equation~\eqref{eq:geodesicraytransformscalar}). As in the case of vector fields in~$\R^n$ there is a gauge for $m$-tensor fields of order $m\geq 1$: if~$h$ is the symmetrized covariant derivative of an $m-1$-tensor field which vanishes on the boundary (or at infinity), then~$h$ is in the kernel of~$\geod_m$. Therefore one can only determine the solenoidal part of the $m$-tensor field from its geodesic ray transform~\cite{IM:integral-geometry-review, PSU-tensor-tomography-progress, SHA-integral-geometry-tensor-fields}; if this can be done we say that~$\geod_m$ is solenoidally injective (or s-injective) on $m$-tensor fields.  

The solenoidal injectivity of~$\geod_m$ has been widely studied and we list only some special cases here: comprehensive treatment can be found in the reviews~\cite{IM:integral-geometry-review, PSU-tensor-tomography-progress}. If $(M, g)$ is a compact simple Riemannian manifold, then the geodesic ray transform is injective on scalar fields and s-injective on 1-forms~\cite{AR-uniqueness-of-one-forms, MU-reconstruction-problem-two-dimensional}.
Solenoidal injectivity is known for tensor fields of any order on two-dimensional compact simple manifolds~\cite{PSU-tensor-tomography-on-simple-surfaces}, on simply connected compact manifolds with strictly convex boundary and non-positive curvature \cite{PS-sharp-stability-nonpositive-curvature, PS-integral-geometry-negative-curvature,SHA-integral-geometry-tensor-fields} and on non-compact Cartan--Hadamard manifolds under certain decay conditions on the tensor fields and on the curvature~\cite{LE-cartan-hadamard, LRS-tensor-tomography-cartan-hadamard}. 
If $n\geq 3$ and $m=0, 1, 2, 4$, then solenoidal injectivity follows from foliation condition by strictly convex hypersurfaces~\cite{dHUZ-inverting-higher-rank-tensors, SUV-inverting-local-tensors, UV-local-geodesic-ray-transform}.
There are also some partial data results for scalar and tensor fields under restrictions on the Riemannian metric~\cite{KRI-support-theorem-analytic-metric, SUV-inverting-local-tensors, UV-local-geodesic-ray-transform}. We also mention that one of the basic general tools in studying solenoidal injectivity of~$\geod_m$ is an energy estimate also known as the Pestov identity~\cite{IM:integral-geometry-review, MU-reconstruction-problem-two-dimensional, PSU-tensor-tomography-progress}.

An interesting generalization of the geodesic ray transform in two dimensions is the mixed ray transform~\cite{deSZ-mixed-ray, SHA-integral-geometry-tensor-fields}
\begin{equation}
\label{eq:mixedraytransform}
\mixed h(\gamma)=\int_{\tau_\gamma^-}^{\tau_\gamma^+}h_{i_1\dotso i_k j_1\dotso j_l}(\gamma (t))(\dot{\gamma}(t)^\perp)^{i_1}\cdots(\dot{\gamma}(t)^\perp)^{i_k}\dot{\gamma}^{j_1}(t)\cdots\dot{\gamma}^{j_l}(t) \der t
\end{equation}
where $\dot{\gamma}(t)^\perp$ denotes the rotation of $\dot{\gamma}(t)$ by 90 degrees counterclockwise and $k+l=m$. The mixed ray transform~$\mixed$ arises in the linearization of the elastic travel time tomography problem~\cite{dESUZ-generic-uniqueness-mixed-ray, deSZ-mixed-ray, SHA-integral-geometry-tensor-fields}.
If $k=0$, then~$\mixed$ reduces to the geodesic ray transform~$\geod_m$. When $l=0$, we have the transverse ray transform~\cite{SHA-integral-geometry-tensor-fields}
\begin{equation}
\label{eq:transverseraytransform}
\geod^\perp_m h(\gamma)=\int_{\tau_\gamma^-}^{\tau_\gamma^+}h_{i_1\dotso i_m}(\gamma (t))(\dot{\gamma}(t)^\perp)^{i_1}\cdots(\dot{\gamma}(t)^\perp)^{i_m} \der t.
\end{equation}
The mixed and transverse ray transforms in equations~\eqref{eq:mixedraytransform} and~\eqref{eq:transverseraytransform} can be extended to higher dimensions $n>2$, but they become tensor-valued transforms~\cite{dESUZ-generic-uniqueness-mixed-ray, SHA-integral-geometry-tensor-fields}. On two-dimensional orientable manifolds one can study the mixed ray transform by reducing it to the geodesic ray transform using rotations~\cite{deSZ-mixed-ray, SHA-integral-geometry-tensor-fields}.

The transverse ray transform was first studied by Braun and Hauck in two-dimensional Euclidean space with applications to flame analysis~\cite{BH-tomographic-reconstruction-vector-fields, NA-mathematical-methods-image-reconstruction, SCHWA-flame-analysis-schlieren, SS-vector-field-overview}. Other applications of the transverse ray transform include diffraction tomography~\cite{LW-diffraction-tomography}, polarization tomography~\cite{SHA-integral-geometry-tensor-fields} and photoelasticity~\cite{HL-applications-to-photoelasticity}. The kernel of the transverse ray transform is known in~$\R^2$ and on higher dimensional manifolds ($n\geq 3$) $\geod^\perp_m$ is even injective under certain conditions~\cite{DS-tomography, NA-mathematical-methods-image-reconstruction, SHA-integral-geometry-tensor-fields}. There are also partial data results for the transverse ray transform~\cite{ABH-support-theorem-transverse-ray, VMS-transverse-partialdata}. For the mixed ray transform some results related to solenoidal injectivity are known in~$\R^2$, on two- and three-dimensional compact simple manifolds, and on manifolds satisfying certain curvature estimates~\cite{dESUZ-generic-uniqueness-mixed-ray, deSZ-mixed-ray, DS-tomography, SHA-integral-geometry-tensor-fields}.

\section{X-ray tomography with partial data: \cite{item:rieszpotential, item:poincareucp, item:covectorpartialdata, item:polynomialucp}}
\label{sec:xraypartialdata}

In the articles~\cite{item:rieszpotential, item:poincareucp, item:covectorpartialdata, item:polynomialucp} we study partial data problems arising in the X-ray tomography of scalar and vector fields. The basic question in such problems is the following: can we say something about the scalar field~$f\colon\R^n\rightarrow\R$ if we know the integrals of~$f$ (the X-ray transform $\xrt_0 f$) on all lines intersecting a given nonempty open set $V\subset\R^n$? We have focused in the uniqueness problem: if $\xrt_0 f=0$ on all lines intersecting~$V$, does it follow that $f=0$? 
In general, the knowledge of the integrals is not enough to determine~$f$ uniquely~\cite{NA-mathematics-computerized-tomography} and therefore one has to put additional assumptions on~$f$. 

We have studied the partial data problem under different assumptions. In the most general case we assume that~$f$ satisfies a constant coefficient partial differential equation in~$V$ in a weak sense. If~$P$ is a polynomial, we let~$P(D)$ be the constant coefficient partial differential operator induced by~$P$, i.e. we consider the partial derivatives~$D$ as variables in~$P$. For example, the polynomial $P(\xi)=\xi_1^2+\dotso+\xi_n^2$ corresponds to the Laplacian $P(D)=-\Delta$.

\begin{figure}[htp]
\centering
\includegraphics[height=6.7cm]{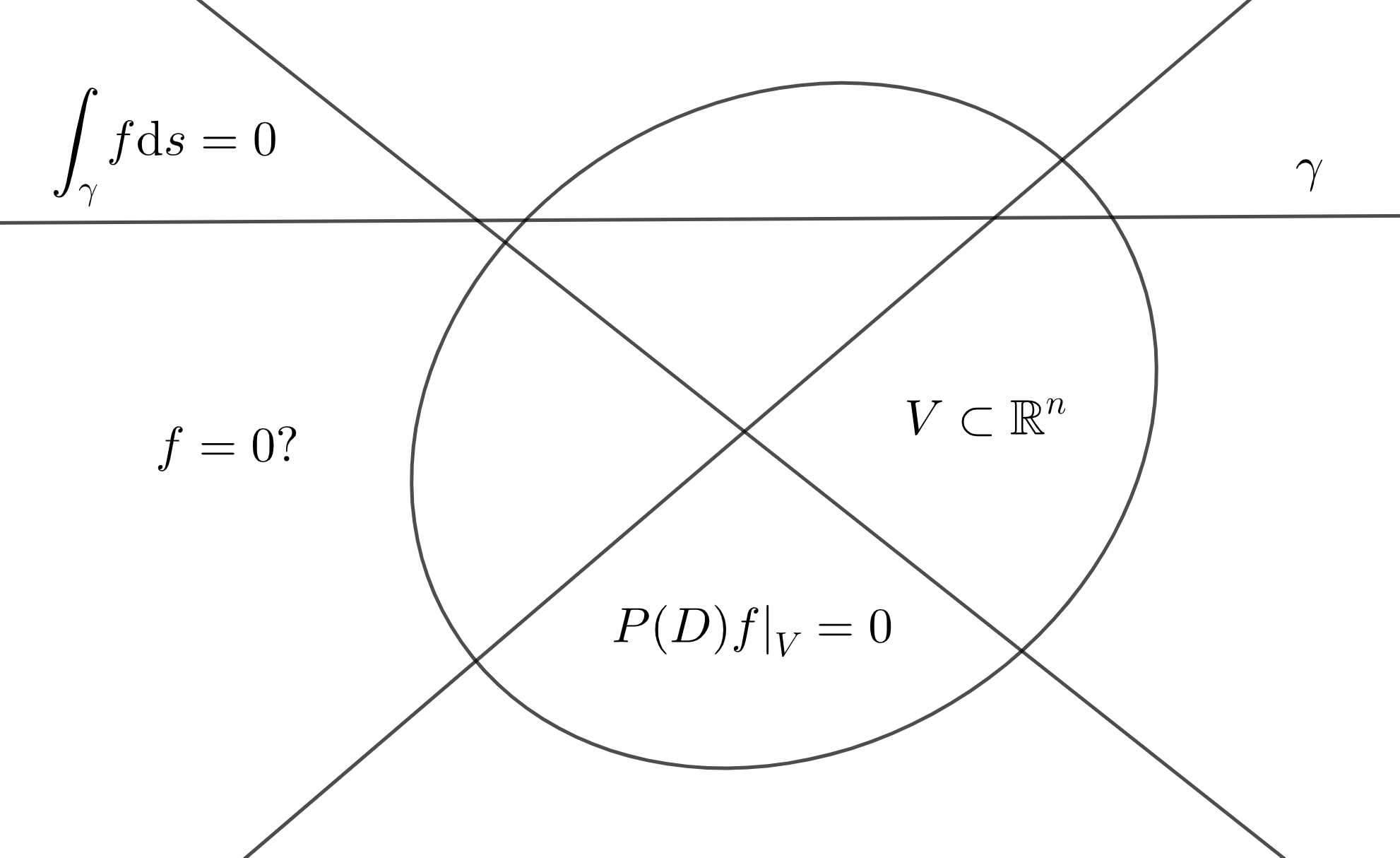}
\caption{The partial data problem for the X-ray transform of scalar fields in its most general form as we have studied. Here $V\subset\R^n$ is a nonempty open set, $P(D)$ is a constant coefficient partial differential operator and~$\gamma$ is a line which intersects~$V$.}
\label{fig:partialdata}
\end{figure}

The main idea of the partial data problem is illustrated in figure~\ref{fig:partialdata}: if $\xrt_0 f=0$ on all lines intersecting~$V$ and $P(D)f|_V=0$ for some constant coefficient partial differential operator~$P(D)$, does it follow that $f=0$ everywhere? Using the linearity of~$\xrt_0$ and~$P(D)$, and the commutativity of distributional derivatives we see that this is indeed a uniqueness problem in the following sense: if~$f_1$ and~$f_2$ are scalar fields such that $P_1(D)f_1|_V=P_2(D)f_2|_V=0$ and $\xrt_0 f_1=\xrt_0 f_2$ on all lines intersecting~$V$, does it follow that $f_1=f_2$ in all of~$\R^n$? The partial data problem can be reduced to a unique continuation problem of the normal operator~$\no_0$ of the X-ray transform: if $\no_0 f|_V=P(D)f|_V=0$, does it follow that $f=0$ everywhere?

More generally, one can replace lines with $d$-planes in the partial data problem of scalar fields. In this way we obtain a partial data problem for the $d$-plane transform~$\dplane$: if $\dplane f=0$ on all $d$-planes intersecting~$V$ and $P(D)f|_V=0$, is it true that $f=0$? 
The partial data problem for vector fields is formulated analogously as in the scalar case. However, for vector fields the problem is naturally formulated in terms of the curl (or the exterior derivative) of the vector field: if $\xrt_1 h=0$ on all lines which intersect~$V$ and $P(D)(\der h)=0$ where~$h$ is a vector field and~$\der h$ its curl, does it follow that $\der h=0$? By the Poincar\'e lemma this is equivalent to that the solenoidal part of~$h$ vanishes~\cite{HO-topological-vector-spaces, MA-poincare-derham-theorems, SHA-integral-geometry-tensor-fields}. As in the case of the X-ray transform of scalar fields, the partial data problems for~$\dplane$ and~$\xrt_1$ can be reduced to the corresponding unique continuation problems of the normal operator~$\nod$ of the $d$-plane transform and the normal operator~$\no_1$ of the X-ray transform of vector fields. We focus on studying the partial data problems from the point of view of the unique continuation of the different normal operators.

In the article~\cite{item:rieszpotential} we study the partial data problem for~$\xrt_0$ under the assumption $f|_V=0$. The main result of the article~\cite{item:rieszpotential} is a unique continuation property of Riesz potentials which correspond to fractional Laplacians with negative exponents. The Riesz potential of a scalar function $f\colon\R^n\to\R$ is defined as the convolution $\riesz f=f\ast\abs{\cdot}^{-\alpha}$ where $\alpha<n$ (see section~\ref{sec:notationxray}). The main theorem is the following: if~$\riesz f$ vanishes to infinite order at some point $x_0\in V$ where the exponent~$\alpha$ satisfies some conditions and $f|_V=0$, then $f=0$. This implies a unique continuation result for the normal operator~$\no_0$: if $\no_0 f|_V=f|_V=0$, then $f=0$. The unique continuation of~$\no_0$ can then be used to prove uniqueness for the partial data problem: if $\xrt_0 f=0$ on all lines intersecting~$V$ and $f|_V=0$, then $f=0$.
We also provide an application of the partial data result to linearized travel time tomography in Euclidean background.

In the article~\cite{item:poincareucp} we study the partial data problem for the $d$-plane transform~$\dplane$ in the case $f|_V=0$. This is a generalization of the problem studied in the article~\cite{item:rieszpotential} where we considered the case $d=1$.
As in the article~\cite{item:rieszpotential}, the partial data problem is studied using the normal operator~$\nod$ of the $d$-plane transform. One of the main results of the article~\cite{item:poincareucp} is a unique continuation property of fractional Laplacians: if $\fraclaplace f|_V=f|_V=0$ where $s\in (-n/2, \infty)\setminus\Z$, then $f=0$. When~$d$ is odd, this implies a unique continuation result for~$\nod$: if $\nod f|_V=f|_V=0$ and~$d$ is odd, then $f=0$. The unique continuation of~$\nod$ then implies uniqueness for the partial data problem: if~$d$ is odd, $\dplane f=0$ on all $d$-planes intersecting~$V$ and $f|_V=0$, then $f=0$. When~$d$ is even and $\dplane f=0$ on all lines which intersect~$V$, we can locally invert the $d$-plane data to obtain that $f|_V=0$.  

The article~\cite{item:covectorpartialdata} considers the partial data problem for~$\xrt_1$ under the assumption $\der h|_V=0$. The approach is similar as in the scalar case, and the main result of the article~\cite{item:covectorpartialdata} is a unique continuation property of~$\no_1$: if~$\no_1 h$ vanishes to infinite order at some point in~$V$ and $\der h|_V=0$, then $\der h=0$. This unique continuation result is proved by reducing it to a unique continuation problem of~$\no_0$ treated in the article~\cite{item:rieszpotential}. The unique continuation of~$\no_1$ can then be used to prove uniqueness for the partial data problem: if $\xrt_1 h=0$ on all lines which intersect~$V$ and $\der h|_V=0$, then $\der h=0$. This is equivalent to that $h=\der\phi$ for some scalar field~$\phi$ by the Poincar\'e lemma, or to that the solenoidal part of~$h$ vanishes. In the article~\cite{item:covectorpartialdata} we also obtain partial data results for the matrix-weighted X-ray transform of vector fields which special case, the Euclidean transverse ray transform, we study in two dimensions.

The article~\cite{item:polynomialucp} is a continuation of the articles~\cite{item:rieszpotential, item:covectorpartialdata} and the integral geometry part of the article~\cite{item:poincareucp}. In particular, we extend the assumptions~$f|_V=0$ and~$\der h|_V=0$ in the partial data problems of scalar and vector fields to the more general cases $P(D)f|_V=0$ and $P(D)(\der h)|_V=0$ where~$P(D)$ is a constant coefficient partial differential operator induced by the polynomial~$P$ as above. The main result of the article~\cite{item:polynomialucp} is a unique continuation property of fractional Laplacians which generalizes the unique continuation result proved in the article~\cite{item:poincareucp}: if $\fraclaplace f|_V=P(D)f|_V=0$ where $s\in (-n/2, \infty)\setminus\Z$ and~$P(D)$ is any constant coefficient partial differential operator, then $f=0$. This unique continuation result directly implies a corresponding unique continuation property for~$\nod$: if~$d$ is odd and $\nod f|_V=P(D)f|_V=0$, then $f=0$. Using reduction to the scalar case one also obtains a unique continuation result for~$\no_1$: if $\no_1 h|_V=P(D)(\der h)|_V=0$, then $\der h=0$. These unique continuation results for~$\nod$ and~$\no_1$ then imply uniqueness for the most general partial data problems we have studied: if $\dplane f=0$ on all $d$-planes intersecting~$V$ where~$d$ is odd and $P(D)f|_V=0$ (or $\xrt_1 h=0$ on all lines intersecting~$V$ and $P(D)(\der h)|_V=0$), then $f=0$ (respectively $\der h=0$).

\subsection{Notation}
\label{sec:notationxray}
Let us first introduce some notation before giving the main theorems. We will follow the notation conventions of the references~\cite{HE:integral-geometry-radon-transforms, ML-strongly-elliptic-systems, MI:distribution-theory2, NA-mathematics-computerized-tomography, SHA-integral-geometry-tensor-fields, SU:microlocal-analysis-integral-geometry, TRE:topological-vector-spaces-distributions}.

We write~$f$ for a scalar function or distribution. The space of tempered distributions is denoted by $\tempered(\R^n)$. We let~$\rapidly(\R^n)\subset\tempered(\R^n)$ be the space of rapidly decreasing distributions. It contains as a subset all compactly supported distributions $\cdistr(\R^n)\subset\rapidly(\R^n)$ and all continuous functions which decrease faster than any polynomial at infinity $C_\infty(\R^n)\subset\rapidly(\R^n)$. The fractional $L^2$-Sobolev space of order $r\in\R$ is defined as
\begin{equation}
H^r(\R^\dimens)=\{f\in\tempered(\R^\dimens): \ifourier(\langle\cdot\rangle^r\hat{f})\in L^2(\R^\dimens)\}
\end{equation}
where $\langle x\rangle=(1+\abs{x}^2)^{1/2}$, $\hat{f}=\fourier(f)$ is the Fourier transform of tempered distributions and~$\ifourier$ is the inverse Fourier transform.
These spaces are nested, i.e. $H^r(\R^n)\hookrightarrow H^t(\R^n)$ continuously when $r\geq t$, and one can identify~$H^{-r}(\R^n)$ with the dual~$(H^r(\R^n))^*$ for every $r\in\R$.
We let $H^{-\infty}(\R^n)=\bigcup_{r\in\R}H^r(\R^n)$ so that $\rapidly(\R^n)\subset H^{-\infty}(\R^n)\subset\tempered(\R^n)$. The fractional Laplacian is defined via Fourier transform
\begin{equation}
\fraclaplace f=\ifourier(\abs{\cdot}^{2s}\hat{f}), \quad s\in (-n/2, \infty)\setminus\Z.
\end{equation}
We have that~$\fraclaplace f$ defines a tempered distribution for $f\in\rapidly(\R^n)$ when $s\in(-n/2, \infty)\setminus\Z$, and for $f\in H^r(\R^n)$ when $s\in(-n/4, \infty)\setminus\Z$.  

The fractional Laplacian has a connection to Riesz potentials. Let $\alpha\in\R$ such that $\alpha<\dimens$. We define the Riesz potential $\riesz\colon\rapidly(\R^\dimens)\rightarrow\tempered(\R^\dimens)$ as $\riesz f=f\ast\kernel$ where the kernel is $\kernel(x)=\abs{x}^{-\alpha}$. If in addition $0<\alpha<\dimens$, then $\riesz=(-\Delta)^{-s}$ up to a constant factor with $s=(\dimens-\alpha)/2$. On the other hand, if $-\dimens/2<s<0$, then we can write $\fraclaplace f=I_{2s+n}f$ up to a constant factor. We say that~$\riesz f$ vanishes to infinite order at a point $x_0\in\R^n$, if~$\riesz f$ is smooth in a neighborhood of~$x_0$ and $\partial^\beta(\riesz f)(x_0)=0$ for all multi-indices $\beta\in\N^n$.

We let~$\pdos$ be the set of all polynomials on~$\R^n$ with complex coefficients excluding the zero polynomial $P\equiv 0$. If $P\in\pdos$ is a polynomial of degree $m\in\N$, then it can be identified with the constant coefficient partial differential operator~$P(D)$ of order $m\in\N$ by writing $P(D)=\sum_{\abs{\alpha}\leq m}a_\alpha D^\alpha$ where $a_\alpha\in\C$, $D^\alpha=D_1^{\alpha_1}\cdots D_n^{\alpha_n}$, $D_j=-i\partial_j$ and $\alpha=(\alpha_1, \dotso, \alpha_n)\in\N^n$ is a multi-index so that $|\alpha|=\alpha_1+\dotso+\alpha_n$. If $V\subset\R^n$ is a nonempty open set, we define the set of admissible functions~$\adm_V$ by setting
\begin{align}
\label{eq:admissible}
\adm_V=\{f\in H^{-\infty}(\R^n): P(D)f|_V=0 \ \text{for some } P\in\pdos\}.
\end{align}
One can see that the set $\adm_V\subset H^{-\infty}(\R^n)$ forms a vector space.

The X-ray transform of scalar fields is denoted by~$\xrt_0$ and it takes a function~$f$ and integrates it over lines. The normal operator is $\no_0=\xrt_0^*\xrt_0$ where~$\xrt_0^*$ is the adjoint of~$\xrt_0$. If~$f$ is a distribution, then~$\xrt_0 f$ and~$\no_0 f$ are defined by duality. More generally, we denote by~$\dplane$ the $d$-plane transform of scalar fields. The transform~$\dplane$ takes a scalar field~$f$ and integrates it over $d$-dimensional planes where $0<d<n$. The normal operator of the $d$-plane transform is $\nod=\dplane^*\dplane$ where~$\dplane^*$ is the adjoint of~$\dplane$. If~$f$ is a distribution, then~$\dplane f$ and~$\nod f$ are defined by using duality.

We denote by~$h$ a vector field or vector-valued distribution. We write $h\in (\cdistr(\R^n))^n$ if $h=(h_1, \dotso, h_n)$ where $h_i\in\cdistr(\R^n)$ for all $i=1, \dotso, n$. The exterior derivative or curl of~$h$ is a matrix whose components are $(\der h)_{ij}=\partial_i h_j-\partial_j h_i$. The X-ray transform of vector fields is denoted by~$\xrt_1$ and it maps a vector field to its line integrals. The normal operator is $\no_1=\xrt_1^*\xrt_1$ where~$\xrt_1^*$ is the adjoint of~$\xrt_1$. If~$h$ is a vector-valued distribution, then both~$\xrt_1 h$ and~$\no_1 h$ are defined by duality.
 
\subsection{Main results}
\label{sec:mainresultsxraypartial}
The following two theorems are the main results of the article~\cite{item:rieszpotential}. The first one is a unique continuation result for Riesz potentials and the second one is a partial data result for the X-ray transform of scalar fields.
\begin{theorem}[{\cite[Theorem 1.1]{item:rieszpotential}}]
\label{thm:rieszpotential}
Let $\alpha=\dimens-1$ or $\alpha\in\R\setminus\Z$ and $\alpha<\dimens$ where $n\geq 2$. Let $f\in\cdistr(\R^{\dimens})$, $V\subset\R^{\dimens}$ any nonempty open set and $x_0\in V$. If $f|_V=0$ and $\riesz f$ vanishes to infinite order at $x_0$, then $f=0$. 
\end{theorem}

\begin{theorem}[{\cite[Theorem 1.2]{item:rieszpotential}}]
\label{thm:partialdatascalar}
Let $V\subset\R^{\dimens}$ be any nonempty open set where $n\geq 2$. If $f\in\cdistr(\R^{\dimens})$ satisfies $f|_V=0$ and $\xrt_0 f$ vanishes on all lines that intersect~$V$, then $f=0$.
\end{theorem}

Theorem~\ref{thm:rieszpotential} is proved by showing that one can obtain all the polynomials in a certain form by taking finite linear combinations of the derivatives of the integral kernel~$\kernel$ in $\riesz f=f\ast\kernel$. The density of polynomials in the space of smooth functions then gives the claim since $f\in\cdistr(\R^n)$ belongs to the dual of that space. We give multiple proofs for theorem~\ref{thm:partialdatascalar}. Two proofs reduce the partial data problem to a unique continuation problem of normal operator: if $\xrt_0 f=0$ on all lines intersecting~$V$, then $\no_0 f|_V=0$. The normal operator~$\no_0$ can be seen as the Riesz potential~$I_{n-1}$ up to a constant factor, or equivalently, as the fractional Laplacian $(-\Delta)^{-1/2}$ up to a constant factor. The partial data result then follows from theorem~\ref{thm:rieszpotential}, or by using the unique continuation of fractional Laplacians which is proved in~\cite{GSU-calderon-problem-fractional-schrodinger}. The third proof works directly at the level of the X-ray transform and is based on angular Fourier series and density of polynomials.

In addition, we provide an application of theorem~\ref{thm:partialdatascalar} to linearized travel time tomography in Euclidean background. In particular, we show how one can use global shear wave splitting data to uniquely determine the difference of the S-wave speeds in weak anisotropy. We also show in the article~\cite{item:rieszpotential} how one can use ``half-local" measurements of travel times to uniquely determine the conformal factor in the linearization; this is a partial data result where we measure travel times of seismic waves in a small open subset of the surface of the Earth, but the waves can emanate from anywhere on the surface.

In the article~\cite{item:poincareucp} we generalize the unique continuation and partial data results proved in~\cite{item:rieszpotential} for scalar fields to $d$-plane transforms. The following two theorems are the main results of the integral geometry part of the article~\cite{item:poincareucp}.

\begin{theorem}[{\cite[Corollary 1]{item:poincareucp}}]
\label{thm:uniquecontinuationdplane}
Let $\dimens \geq 2$ and let $f$ belong to either~$\cdistr(\R^\dimens)$ or~$C_{\infty}(\R^\dimens)$. Let $d \in \N$ be odd such that $0<d<\dimens$. If $\nod f|_V=0$ and $f|_V=0$ for some nonempty open set $V\subset\R^\dimens$, then $f=0$.
\end{theorem}

\begin{theorem}[{\cite[Corollary 2]{item:poincareucp}}]
\label{thm:partialdatadplane}
Let $\dimens \geq 2$, $V\subset\R^\dimens$ a nonempty open set and $f\in C_{\infty}(\R^\dimens)$ or $f\in\cdistr(\R^\dimens)$. Let $d\in\N$ be odd such that $0<d<\dimens$. If $f|_V=0$ and $\dplane f=0$ for all $d$-planes intersecting $V$, then $f=0$. 
\end{theorem}

Theorem~\ref{thm:uniquecontinuationdplane} is proved by using a unique continuation property of fractional Laplacians which is proved in the same article~\cite{item:poincareucp} (see theorem~\ref{thm:ucplaplacian}). Unique continuation of fractional Laplacians can be used since the normal operator~$\nod$ of the $d$-plane transform corresponds to the fractional Laplacian~$(-\Delta)^{-d/2}$ up to a constant factor. The unique continuation of~$\nod$ is then used to prove theorem~\ref{thm:partialdatadplane}. For this reason we have to assume that~$d$ is odd: theorem~\ref{thm:uniquecontinuationdplane} does not hold if~$d$ is even since in that case~$\nod$ is the inverse of a local operator. However, if~$d$ is even, then the partial data problem for the $d$-plane transform is locally uniquely solvable: if $\dplane f=0$ on all lines intersecting~$V$ and~$d$ is even, then $f|_V=0$.

In the article~\cite{item:covectorpartialdata} we generalize the above partial data results to vector fields. The following two main theorems of the article~\cite{item:covectorpartialdata} are similar to theorems~\ref{thm:rieszpotential} and~\ref{thm:partialdatascalar}. 

\begin{theorem}[{\cite[Theorem 1.1]{item:covectorpartialdata}}]
\label{thm:uniquecontinuationvector}
Let $h\in (\cdistr(\R^\dimens))^\dimens$ and $V\subset\R^\dimens$ some nonempty open set where $n\geq 2$.
If $\der h|_V=0$ and~$\no_1 h$ vanishes to infinite order at $x_0\in V$, then $h=\der\phi$ for some $\phi\in\cdistr(\R^\dimens)$.
\end{theorem}

\begin{theorem}[{\cite[Theorem 1.2]{item:covectorpartialdata}}]
\label{thm:partialdatavector}
Let $h\in (\cdistr(\R^\dimens))^\dimens$ and $V\subset\R^\dimens$ some nonempty open set where $n\geq 2$. Assume that $\der h|_V=0$.
Then $\xrt_1 h$ vanishes on all lines intersecting~$V$ if and only if $h=\der\phi$ for some $\phi\in\cdistr (\R^\dimens)$. 
\end{theorem}

Instead of assuming that~$\no_1 h$ vanishes to infinite order in theorem~\ref{thm:uniquecontinuationvector} we could require that~$\der(\no_1 h)$ vanishes componentwise to infinite order at some point $x_0\in V$. This weaker condition implies the claim since theorem~\ref{thm:uniquecontinuationvector} is proved by using theorem~\ref{thm:rieszpotential} and the fact that $\der(\no_1 h)=\no_0(\der h)$ holds componentwise up to a constant factor. We provide two alternative proofs for theorem~\ref{thm:partialdatavector}. The first proof directly uses the unique continuation of the normal operator~$\no_1$ in theorem~\ref{thm:uniquecontinuationvector}. 
The second proof is based on Stokes' theorem and theorem~\ref{thm:partialdatascalar}. Both proofs use the same idea: from the assumptions we deduce that $\der h=0$ and the Poincar\'e lemma implies that $h=\der\phi$ for some scalar field~$\phi$. 

In the article~\cite{item:covectorpartialdata} we also study the matrix-weighted X-ray transform of vector fields $\xrt_A=\xrt_1\circ A$ where~$A$ is a smooth invertible matrix field. Similar results as in theorems~\ref{thm:uniquecontinuationvector} and~\ref{thm:partialdatavector} are obtained for the transform~$\xrt_A$. As a special case of the transform~$\xrt_A$ we obtain results for the Euclidean transverse ray transform in two dimensions.

In the article~\cite{item:polynomialucp} we generalize the partial data and unique continuation results obtained in the articles~\cite{item:rieszpotential, item:poincareucp, item:covectorpartialdata}. The partial data results are proved by using the following unique continuation property of fractional Laplacians which is a generalization of the unique continuation result we proved in the article~\cite{item:poincareucp}.

\begin{theorem}[{\cite[Theorem 1.1]{item:polynomialucp}}]
\label{thm:ucppolynomial}
Let $n\geq 1$, $s\in(-n/4, \infty)\setminus\Z$ and $f\in\adm_V$ where $V\subset\R^n$ is some nonempty open set. If $\fraclaplace f|_V=0$, then $f=0$. If $f\in\rapidly(\R^n)\cap\adm_V$, then the claim holds for $s\in(-n/2, \infty)\setminus\Z$.
\end{theorem}

The condition $f\in\adm_V$ means that $f\in H^r(\R^n)$ for some $r\in\R$ and $P(D)f|_V=0$ for some constant coefficient partial differential operator~$P(D)$.
Theorem~\ref{thm:ucppolynomial} is proved by using the unique continuation result of fractional Laplacians proved in~\cite{item:poincareucp} (see theorem~\ref{thm:ucplaplacian}) for the scalar field~$P(D)f$. The assumptions and locality of~$P(D)$ imply the conditions $P(D)f|_V=\fraclaplace (P(D)f)|_V=0$ and hence~$f$ has to satisfy the global partial differential equation $P(D)f=0$ which has only trivial solutions in the class of admissible functions~$\adm_V$.

As before, the unique continuation of fractional Laplacians in theorem~\ref{thm:ucppolynomial} can be used to prove partial data results for scalar and vector fields. The following two theorems of the article~\cite{item:polynomialucp} are generalizations of theorems~\ref{thm:partialdatascalar} and~\ref{thm:partialdatavector}.

\begin{theorem}[{\cite[Theorem 1.4]{item:polynomialucp}}]
\label{thm:xrayscalarpolynomial}
Let $n\geq 2$ and $f\in\cdistr(\R^n)\cap\adm_V$ or $f\in C_\infty(\R^n)\cap\adm_V$ where $V\subset\R^n$ is some nonempty open set. If $\xrt_0 f=0$ on all lines intersecting~$V$, then $f=0$.
\end{theorem}

\begin{theorem}[{\cite[Theorem 1.7]{item:polynomialucp}}]
\label{thm:xrayvectorpolynomial}
Let $n\geq 2$ and $h\in (\cdistr(\R^n))^n$ such that $(\der h)_{ij}\in\adm_V$ for all $i, j=1, \dotso, n$ where $V\subset\R^n$ is some nonempty open set. If $\xrt_1 h=0$ on all lines intersecting~$V$, then $\der h=0$. Especially, $h=\der\phi$ for some $\phi\in\cdistr(\R^n)$.
\end{theorem}

Theorems~\ref{thm:xrayscalarpolynomial} and~\ref{thm:xrayvectorpolynomial} are proved in the following way.
Theorem~\ref{thm:ucppolynomial} implies a corresponding unique continuation result for the normal operator~$\no_0$. The unique continuation of~$\no_0$ is then used to prove the partial data result in theorem~\ref{thm:xrayscalarpolynomial}. Further, using again the fact that $\der (\no_1 h)=\no_0(\der h)$ holds componentwise up to a constant factor we can prove a unique continuation property for~$\no_1$, which in turn implies the partial data result in theorem~\ref{thm:xrayvectorpolynomial}. When~$d$ is odd, theorem~\ref{thm:ucppolynomial} implies a corresponding unique continuation result for the normal operator~$\nod$, which in turn implies a similar partial data result as in theorem~\ref{thm:xrayscalarpolynomial} for the $d$-plane transform~$\dplane$.

\section{Higher order fractional Calder\'on problems: \cite{item:poincareucp, item:gunther}}
\label{sec:fractionalcalderon}

In the articles~\cite{item:poincareucp, item:gunther} we study uniqueness for higher order fractional Calder\'on problems. Let $\Omega\subset\R^n$ be an open bounded set, $\Omega_e=\R^n\setminus\overline{\Omega}$ its exterior and $s\in (0, \infty)\setminus\Z$. We consider the Calder\'on problem for the fractional Schr\"odinger equation 
\begin{align}\label{eq:schrodinger} \left\{\begin{array}{rl}
        (\fraclaplace+q)u&=0  \;\;\text{in}\;  \Omega\\
        u|_{\Omega_e} &=f
        \end{array}\right.\end{align}
and for the more general equation involving lower order local perturbations of the fractional Laplacian
\begin{align}\label{eq:perturbations} \left\{\begin{array}{rl}
        (\fraclaplace+P(x, D))u&=0  \;\;\text{in}\;  \Omega\\
        u|_{\Omega_e} &=f
        \end{array}\right.\end{align}
where~$P(x, D)$ is a variable coefficient partial differential operator of order $m\in\N$. We can write~$P(x, D)$ as
\begin{equation}
P(x, D)=\sum_{\abs{\alpha}\leq m}a_{\alpha}(x)D^{\alpha}
\end{equation}
where the coefficients $a_\alpha=a_\alpha(x)$ are functions in~$\Omega$ (or more generally Sobolev multipliers in~$\R^n$).
We assume that $m<2s$ so that~$P(x, D)$ can be considered as a lower order perturbation to~$\fraclaplace$. We see that equation~\eqref{eq:schrodinger} is a special case of equation~\eqref{eq:perturbations} and the potential~$q$ can be treated as a zeroth order perturbation to~$\fraclaplace$.

The inverse problem for equations~\eqref{eq:schrodinger} and~\eqref{eq:perturbations} is illustrated in figure~\ref{fig:calderon}. Formally, we put some ``non-local voltage" in the open set $W_1\subset\Omega_e$ and measure ``non-local currents" in the open set $W_2\subset\Omega_e$. More precisely, the fractional Calder\'on problem for equation~\eqref{eq:perturbations} is formulated as follows: if the DN maps~$\Lambda_{P_1}$ and~$\Lambda_{P_2}$ agree in~$W_2$ for all exterior values $f\in C_c^\infty(W_1)$, does it follow that the partial differential operators~$P_1$ and~$P_2$ are equal in~$\Omega$?
This problem was first introduced by Ghosh, Salo and Uhlmann in their seminal work~\cite{GSU-calderon-problem-fractional-schrodinger} where the authors studied equation~\eqref{eq:schrodinger} in the case $s\in (0, 1)$. We can think the inverse problem as a partial data problem: instead of having data in the full exterior~$\Omega_e$ we only have information or make measurements in the (possibly small) open subsets $W_1, W_2\subset\Omega_e$.

\begin{figure}[htp]
\centering
\includegraphics[height=6.5cm]{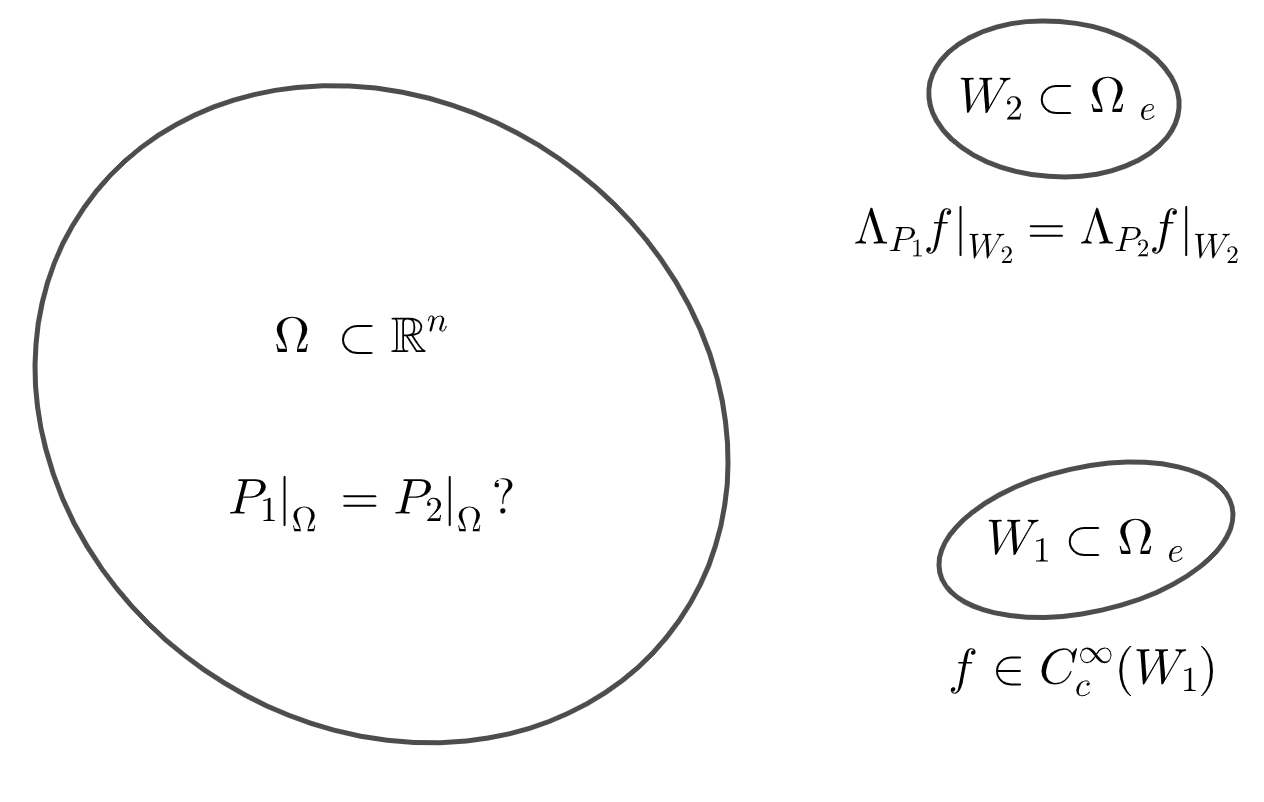}
\caption{The fractional Calder\'on problem in its most general form as we have studied. Here $\Omega\subset\R^n$ is a bounded open set and $\Omega_e=\R^n\setminus\overline{\Omega}$ its exterior. The ``measurements" are done in the (possibly disjoint) subsets $W_1, W_2\subset\Omega_e$ of the exterior.}
\label{fig:calderon}
\end{figure}

The basic tools in proving uniqueness for fractional Calder\'on problems are the unique continuation property and fractional Poincar\'e inequality for the fractional Laplacian~$\fraclaplace$. We already saw the importance of unique continuation of fractional Laplacians in partial data problems of scalar and vector fields. In fractional Calder\'on problems the unique continuation of~$\fraclaplace$ implies Runge approximation: one can approximate functions in certain Sobolev spaces arbitrarily well by solutions of the fractional equation under study (see section~\ref{sec:mainresultsfractional}). The fractional Poincar\'e inequality is a norm estimate involving the $L^2$-norms of a function and its fractional Laplacian, and it is an important inequality in proving well-posedness for the forward problem (the coercivity of the bilinear form). These basic tools (unique continuation and Poincar\'e inequality) were proved in~\cite{GSU-calderon-problem-fractional-schrodinger} in the case $s\in (0, 1)$.  

In the article~\cite{item:poincareucp} we study the higher order fractional Calder\'on problem for equation~\eqref{eq:schrodinger} when $s\in (0, \infty)\setminus\Z$. We prove higher order unique continuation result for fractional Laplacians: if $\fraclaplace u|_V=u|_V=0$ for some nonempty open set $V\subset\R^n$ and $s\in (-n/2, \infty)\setminus\Z$, then $u=0$. This generalizes the result proved in~\cite{GSU-calderon-problem-fractional-schrodinger}. We also prove higher order fractional Poincar\'e inequality for $s\in (0, \infty)\setminus\Z$ which says that the $L^2$-norm of~$u$ can be bounded from above by the $L^2$-norm of~$\fraclaplace u$. We provide five possible proofs for the Poincar\'e inequality and some of the proofs also give information about the constant in the inequality. Using the Poincar\'e inequality we prove well-posedness of the forward problem, and unique continuation of~$\fraclaplace$ implies Runge approximation for equation~\eqref{eq:schrodinger}. Using Runge approximation and the so-called Alessandrini identity (see section~\ref{sec:mainresultsfractional}) for suitable test functions we prove uniqueness for the inverse problem: if $W_1, W_2\subset\Omega_e$ are some open sets such that the DN maps satisfy $\Lambda_{q_1}f=\Lambda_{q_2}f$ in~$W_2$ for all exterior values~$f\in C_c^\infty(W_1)$, then $q_1=q_2$ in~$\Omega$. This is done for certain singular potentials~$q$ which can be viewed as Sobolev multipliers. We also study the magnetic counterpart of equation~\eqref{eq:schrodinger} (the higher order fractional magnetic Schr\"odinger equation) in the article~\cite{item:poincareucp} and prove uniqueness (up to a gauge) under certain assumptions on the electric and magnetic potentials, generalizing the results in~\cite{CO-magnetic-fractional-schrodinger} to higher order cases.

The article~\cite{item:gunther} is a continuation of the article~\cite{item:poincareucp} in the sense that we replace the potential~$q$ in equation~\eqref{eq:schrodinger} with a general lower order local perturbation $P(x, D)$. In the article~\cite{item:gunther} we study the fractional Calder\'on problem for equation~\eqref{eq:perturbations} when $s\in (0, \infty)\setminus\Z$ and $m<2s$.
We consider two different classes of coefficients~$a_\alpha$ of the partial differential operator~$P(x, D)$: coefficients which belong to certain $L^\infty$-Bessel potential spaces, and coefficients which are certain Sobolev multipliers. The same tools that we develop in the article~\cite{item:poincareucp}, i.e. the higher order unique continuation property and fractional Poincar\'e inequality for~$\fraclaplace$, are applicable in proving uniqueness in the article~\cite{item:gunther}. In addition to the Poincar\'e inequality we also need the Kato--Ponce inequality in proving well-posedness of the forward problem. The Kato--Ponce inequality is a fractional Leibnitz rule in terms of $L^p$-norms~\cite{GO14, GK-schrodinger-operators, KAPO-commutator-estimates-kato-ponce}. As in the article~\cite{item:poincareucp}, the unique continuation of~$\fraclaplace$ implies Runge approximation for equation~\eqref{eq:perturbations}. Using the Runge approximation and the corresponding Alessandrini identity for suitable test functions we prove uniqueness for the inverse problem: if $W_1, W_2\subset\Omega_e$ are some open sets such that the DN maps satisfy $\Lambda_{P_1}f=\Lambda_{P_2}f$ in~$W_2$ for all exterior values~$f\in C_c^\infty(W_1)$, then $P_1=P_2$ in~$\Omega$. This uniqueness result is shown for both classes of coefficients~$a_\alpha$, i.e. coefficients with bounded fractional derivatives and coefficients which are Sobolev multipliers.
\subsection{Notation}
\label{sec:notationfractional}
We use the same notation that we introduced in section~\ref{sec:notationxray}, but we also introduce some additional notation. We follow the notation conventions of the references~\cite{BL-interpolation-spaces, CWHM-sobolev-spaces-on-non-lipchtiz-domains, MS-theory-of-sobolev-multipliers, ML-strongly-elliptic-systems, MI:distribution-theory2, TRI-interpolation-function-spaces}.

If $1\leq p\leq\infty$, we define the fractional $L^p$-Bessel potential space of order $r\in\R$ as
\begin{equation}
H^{r, p}(\R^\dimens)=\{u\in\tempered(\R^\dimens): \ifourier(\langle\cdot\rangle^r\hat{u})\in L^p(\R^\dimens)\}
\end{equation}
and we equip it with the norm
\begin{equation}
\aabs{u}_{H^{r,p}(\R^\dimens)}=\aabs{\ifourier(\langle\cdot\rangle^r\hat{u})}_{L^p(\R^\dimens)}.
\end{equation}
These spaces are nested, i.e. $H^{r, p}(\R^\dimens)\hookrightarrow H^{t, p}(\R^\dimens)$ continuously when $r\geq t$. We see that $H^{0, p}(\R^n)=L^p(\R^n)$.
If $\Omega\subset\R^n$ is an open set, then we define the spaces~$H^{r, p}(\Omega)$ as restrictions
\begin{equation}
H^{r, p}(\Omega)=\{u|_\Omega: u\in H^{r, p}(\R^\dimens)\}
\end{equation}
and we use the quotient norm
\begin{equation}
\aabs{w}_{H^{r, p}(\Omega)}=\inf \{\aabs{u}_{H^{r, p}(\R^\dimens)}: u\in H^{r, p}(\R^\dimens) \ \text{such that }  u|_\Omega=w\}.
\end{equation}
It follows that the inclusions $H^{r, p}(\Omega)\hookrightarrow H^{t, p}(\Omega)$ are continuous when $r\geq t$. The spaces~$H^{r, p}(\Omega)$ are not to be confused with the Sobolev-Slobodeckij spaces~$W^{r, p}(\Omega)$ which are defined by using weak derivatives of $L^p$-functions and which in general are different from the Bessel potential spaces we have introduced~\cite{DINEPV-hitchhiker-sobolev}. If $r\geq 0$ and $p=2$, then $H^{r, 2}(\R^n)=W^{r, 2}(\R^n)$ and  $H^{r, 2}(\Omega)=W^{r, 2}(\Omega)$ when~$\Omega$ is a Lipschitz domain.

The following spaces are special cases of the above Bessel potential spaces
\begin{align}
H_F^{r, p}(\R^\dimens) &=\{u\in H^{r, p}(\R^\dimens): \spt(u)\subset F\} \\
\widetilde{H}^{r, p}(\Omega)&= \ \text{closure of} \ C_c^{\infty}(\Omega) \ \text{with respect to the norm } \aabs{\cdot}_{H^{r, p}(\R^n)} \\
H_0^{r, p}(\Omega)&= \ \text{closure of} \ C_c^{\infty}(\Omega) \ \text{with respect to the norm } \aabs{\cdot}_{H^{r, p}(\Omega)}
\end{align}
where $F\subset\R^n$ is some closed set. Observe that $\widetilde{H}^{r, p}(\Omega)\subset H^{r, p}(\R^n)$ and $H_0^{r, p}(\Omega)\subset H^{r, p}(\Omega)$.
One also sees that $\widetilde{H}^{r, p}(\Omega)\subset H_0^{r, p}(\Omega)$ and $\widetilde{H}^{r, p}(\Omega)\subset H^{r, p}_{\overline{\Omega}}(\R^\dimens)$.
When $p=2$, we simply write $H^{r, 2}(\R^n)=H^r(\R^n)$, $H^{r, 2}(\Omega)=H^r(\Omega)$ and so on. It follows that $(\widetilde{H}^r(\Omega))^*=H^{-r}(\Omega)$ and $(H^r(\Omega))^*=\widetilde{H}^{-r}(\Omega)$ for any open set $\Omega\subset\R^n$ and $r\in\R$. If in addition~$\Omega$ is a Lipschitz domain and $r\geq 0$ such that $r\notin \{\frac{1}{2}, \frac{3}{2}, \frac{5}{2} \dotso\}$, then $\widetilde{H}^r(\Omega)=H^r_0(\Omega)$. 

We define the space of Sobolev multipliers $M(H^r\rightarrow H^t)\subset\distr(\R^n)$ by saying that the distribution $f\in\distr(\R^n)$ belongs to $M(H^r\rightarrow H^t)$ if the multiplier norm 
\begin{equation}
\|f\|_{r,t} = \sup \{\abs{\ip{f}{u v}}: u,v \in C_c^\infty(\mathbb R^n), \ \aabs{u}_{H^r(\R^n)} = \aabs{v}_{H^{-t}(\R^n)} =1 \}
\end{equation}
is finite. We let $M_0(H^r\rightarrow H^t)$ be the closure of $C^\infty_c(\R^n)$ in $M(H^r\rightarrow H^t)$ with respect to the norm $\aabs{\cdot}_{r, t}$. The elements of the space~$M(H^r\rightarrow H^t)$ are called Sobolev multipliers since each $f \in M(H^r\rightarrow H^t)$ induces a map $m_f\colon H^r(\R^n) \rightarrow H^t(\R^n)$ defined as
\begin{align}
    \langle m_f(u),v \rangle = \langle f,uv\rangle
\end{align}
for all $u\in H^r(\R^n)$ and $\ v \in H^{-t}(\R^n)$.
As a special case of multipliers we write $Z^{-s}(\R^n)=M(H^s\rightarrow H^{-s})$ and $Z^{-s}_0(\R^n)=M_0(H^s\rightarrow H^{-s})$ whose elements we also call singular potentials. 

We say that~$0$ is not a Dirichlet eigenvalue of the operator $\fraclaplace+q$, if the following condition holds:
\begin{equation}
\label{eq:zeronotdirichleteigenvalue}
\text{If} \ u\in H^s(\R^\dimens) \ \text{solves} \ (\fraclaplace+q)u=0 \ \text{in} \ \Omega \ \text{and} \ u|_{\Omega_e}=0, \ \text{then} \ u=0.
\end{equation}
Analogously, we say that~$0$ is not a Dirichlet eigenvalue of the operator $\fraclaplace+P(x, D)$ if condition~\eqref{eq:zeronotdirichleteigenvalue} holds when~$q$ is replaced with the partial differential operator~$P(x, D)$. When the forward problem for equation~\eqref{eq:schrodinger} is well-posed, we can define the DN map~$\Lambda_q\colon H^s(\Omega_e)\to (H^s(\Omega_e))^*$ as $\ip{\Lambda_q f_1}{f_2}=B_q(u_{f_1}, f_2)$ where $B_q(\cdot, \cdot)$ is the bilinear form associated to equation~\eqref{eq:schrodinger} and~$u_{f_1}$ is the unique solution to equation~\eqref{eq:schrodinger} with exterior value $u|_{\Omega_e}=f_1$. The DN map~$\Lambda_P$ for equation~\eqref{eq:perturbations} is defined similarly.

\subsection{Main results}
\label{sec:mainresultsfractional}
One of the main theorems of the article~\cite{item:poincareucp} is the following unique continuation property of fractional Laplacians.
\begin{theorem}[{\cite[Theorem 1.1]{item:poincareucp}}]
\label{thm:ucplaplacian}
Let $\dimens \geq 1$, $s\in (-n/4,\infty)\setminus \Z$ and $u\in H^{r}(\R^\dimens)$ where $r\in\R$. If $(-\Delta)^s u|_V=0$ and $u|_V=0$ for some nonempty open set $V\subset\R^\dimens$, then $u=0$. The claim holds also for $s\in (-n/2, -n/4]\setminus\Z$ if $u\in H^{r, 1}(\R^\dimens)$ or $u\in\rapidly(\R^\dimens)$.
\end{theorem}

Theorem~\ref{thm:ucplaplacian} is proved by reducing the claim to the case $s\in (0, 1)$ and using the unique continuation result proved in~\cite{GSU-calderon-problem-fractional-schrodinger}. The reduction can be done by using the simple relation $(-\Delta)^k\fraclaplace=\fraclaplace (-\Delta)^k=(-\Delta)^{s+k}$ when $k\in\N$ and $s\in(-n/2, \infty)\setminus\Z$. The assumptions on~$s$ in theorem~\ref{thm:ucplaplacian} are put so that~$\fraclaplace$ is a non-local operator and~$\fraclaplace u$ is well-defined as a tempered distribution. We also prove in the article~\cite{item:poincareucp} many other versions of the unique continuation of~$\fraclaplace$ in different Sobolev spaces, including homogeneous Sobolev spaces and certain Bessel potential spaces.

The next theorem of the article~\cite{item:poincareucp} is called the (fractional) Poincar\'e inequality. It has an essential role in proving well-posedness for the forward problems of equations~\eqref{eq:schrodinger} and~\eqref{eq:perturbations}. 

\begin{theorem}[{\cite[Theorem 1.2]{item:poincareucp}}]
\label{thm:poincare}
Let $\dimens \geq 1$, $s\geq t\geq 0$, $K\subset\R^\dimens$ a compact set and $u\in H_K^s(\R^\dimens)$. There exists a constant $c=c(n, K, s)> 0$ such that
\begin{equation}
\label{eq:poincareinequality}
\aabs{(-\Delta)^{t/2}u}_{L^2(\R^\dimens)}\leq c\aabs{(-\Delta)^{s/2}u}_{L^2(\R^\dimens)}.
\end{equation}
\end{theorem}

In well-posedness we only need the cases $t=0$ and $s\in (0, \infty)\setminus\Z$ of theorem~\ref{thm:poincare}. Note that the inequality~\eqref{eq:poincareinequality} holds for all exponents $s\geq t\geq 0$, not just fractional ones. The interpretation of theorem~\ref{thm:poincare} is that the norms of lower order derivatives of~$u$ are bounded from above by the norms of higher order derivatives of~$u$. When $t=0$ and $s=1$, then the inequality~\eqref{eq:poincareinequality} reduces to the classical Poincar\'e inequality. 

We provide five different proofs for theorem~\ref{thm:poincare}. Two of the simplest proofs are based on Fourier analysis: the first uses splitting of frequencies on the Fourier side and the second uses uncertainty inequalities proved in~\cite{FS-the-uncertainty-principle}. Two other proofs are based on a reduction argument similar to what we did in proving the unique continuation of higher order fractional Laplacians. The fifth proof considers the case $s\geq 1$ and it uses interpolation in homogeneous Sobolev spaces and the classical Poincar\'e inequality. This proof also gives an explicit constant for the inequality: if $s\geq 1$ and $u\in\widetilde{H}^s(\Omega)$, then in theorem~\ref{thm:poincare} we can take $c=C^{s-t}$  where~$C$ is the classical Poincar\'e constant. This is expected since on the left-hand side of equation~\eqref{eq:poincareinequality} we take~$t$ derivatives and on the right-hand side we take~$s$ derivatives.

The next theorem of the article~\cite{item:poincareucp} gives uniqueness for the higher order fractional Schr\"odinger equation with a singular potential.

\begin{theorem}[{\cite[Theorem 1.3]{item:poincareucp}}]
\label{thm:higherorderuniqueness}
Let $\dimens \geq 1$, $\Omega\subset\R^\dimens$ a bounded open set, $s\in (0, \infty)\setminus\Z$, and $q_1, q_2\in Z_0^{-s}(\R^\dimens)$ such that~$0$ is not a Dirichlet eigenvalue of the operators $\fraclaplace+q_i$. Let $W_1, W_2\subset\R^n\setminus\overline{\Omega}$ be open sets. If the DN maps for the equations $\fraclaplace u+m_{q_i}(u)=0$ in $\Omega$ satisfy $\Lambda_{q_1}f|_{W_2}=\Lambda_{q_2}f|_{W_2}$ for all $f\in C_c^{\infty}(W_1)$, then $q_1|_{\Omega}=q_2|_{\Omega}$.
\end{theorem}

The proof of theorem~\ref{thm:higherorderuniqueness} follows from Runge approximation for equation~\eqref{eq:schrodinger} and choosing suitable test functions in the Alessandrini identity. The Runge approximation says that we can approximate functions in $\widetilde{H}^s(\Omega)$ arbitrarily well by solutions of the fractional Schr\"odinger equation~\eqref{eq:schrodinger}, and it can be proved by using the unique continuation of~$\fraclaplace$ in theorem~\ref{thm:ucplaplacian} and the well-posedness of the forward problem. The Alessandrini identity is an integral identity showing how the DN maps~$\Lambda_{q_i}$ and the corresponding potentials~$q_i$ are related in terms of exterior values~$f$ and solutions~$u_f$ of equation~\eqref{eq:schrodinger}.

The article~\cite{item:gunther} generalizes the higher order fractional Schr\"odinger equation studied in the article~\cite{item:poincareucp} to include more general lower order local perturbations. The following two theorems are the main results of the article~\cite{item:gunther}.

\begin{theorem}[{\cite[Theorem 1.1]{item:gunther}}]
\label{thm:perturbations1}
Let $\Omega \subset \mathbb R^n$ be a bounded open set where $n\geq 1$. Let $s \in (0, \infty) \setminus \mathbb Z$ and $m\in \mathbb N$ be such that $2s > m$. Let 
\begin{equation}
P_j = \sum_{|\alpha|\leq m} a_{j,\alpha} D^\alpha,\quad j = 1, 2,
\end{equation} 
be partial differential operators of order $m$ where $a_{j,\alpha} \in M_0(H^{s-|\alpha|}\rightarrow H^{-s})$ such that~$0$ is not a Dirichlet eigenvalue of the operators $\fraclaplace+P_j$. Given any two open sets $W_1, W_2 \subset\R^n\setminus\overline{\Omega}$, suppose that the DN maps $\Lambda_{P_j}$ for the equations
$((-\Delta)^s + P_j)u = 0$ in $\Omega$ satisfy 
\begin{equation}
\Lambda_{P_1}f|_{W_2} = \Lambda_{P_2}f|_{W_2}
\end{equation}
for all $f \in C^\infty_c(W_1)$. Then $P_1|_{\Omega} = P_2|_{\Omega}$.
\end{theorem}

\begin{theorem}[{\cite[Theorem 1.2]{item:gunther}}]
\label{thm:perturbations2} Let $\Omega \subset \mathbb R^n$ be a bounded Lipschitz domain where $n\geq 1$. Let $s \in (0, \infty) \setminus \mathbb Z$ and $m\in \mathbb N$ be such that $2s > m$. Let 
\begin{equation}
P_j (x, D) = \sum_{|\alpha|\leq m} a_{j,\alpha}(x) D^\alpha,\quad j = 1, 2,
\end{equation}
be partial differential operators of order $m$ with coefficients $a_{j,\alpha} \in H^{r_\alpha,\infty}(\Omega)$ where
\begin{align}\label{r-alpha-conditions}
    r_\alpha = \Bigg\{ \begin{matrix*} 
    0 & \mbox{if} & |\alpha|-s<0, \\ 
    |\alpha|-s+\delta & \mbox{if} & |\alpha|-s \in \{ 1/2, 3/2, ... \}, \\
    |\alpha|-s & \mbox{if} & \mbox{otherwise} \\
    \end{matrix*} \Bigg. 
\end{align}
for any fixed $\delta >0$ and assume that~$0$ is not a Dirichlet eigenvalue of the operators $\fraclaplace+ P_j(x, D)$. Given any two open sets $W_1, W_2 \subset\R^n\setminus\overline{\Omega}$, suppose that the DN maps $\Lambda_{P_j}$ for the equations
$((-\Delta)^s + P_j (x, D))u = 0$ in $\Omega$ satisfy 
\begin{equation}
\Lambda_{P_1}f|_{W_2} = \Lambda_{P_2}f|_{W_2}
\end{equation}
for all $f \in C^\infty_c(W_1)$. Then $P_1(x, D) = P_2(x, D)$.
\end{theorem}

It is not known whether the spaces $M_0(H^{s-|\alpha|}\rightarrow H^{-s})$ and~$H^{r_\alpha,\infty}(\Omega)$ are contained in each other. If this is not the case, then theorems~\ref{thm:perturbations1} and~\ref{thm:perturbations2} are distinct and neither claim implies the other. In theorem~\ref{thm:perturbations1} we consider multipliers which can be approximated in the multiplier norm by smooth compactly supported functions and for this reason we do not need to assume anything about the boundary of~$\Omega$. In theorem~\ref{thm:perturbations2} we have put some conditions on~$\partial\Omega$ and for the exponent~$r_\alpha$ which are needed in proving well-posedness in the case of coefficients with bounded fractional derivatives. The assumptions that~$0$ is not a Dirichlet eigenvalue and $2s>m$ (i.e. we consider perturbations to~$\fraclaplace$) are also crucial in both theorems when proving well-posedness of the forward problem. 
It follows that $M(H^{s-|\alpha|}\rightarrow H^{-s}) = \{0\}$ if $s-\abs{\alpha} < -s$. Partly because of this reason theorem~\ref{thm:perturbations1} is formulated only for $2s>m$ since the multiplier coefficients for higher order derivatives are zero, i.e. $a_\alpha = 0$ for all $\abs{\alpha} > 2s$. 

Theorems~\ref{thm:perturbations1} and~\ref{thm:perturbations2} are proved in the same way, even though the exact details are a little bit different. The proofs follow the same ideas as in the article~\cite{item:poincareucp} where we proved uniqueness for zeroth order perturbations, and we see that theorem~\ref{thm:higherorderuniqueness} is in fact a special case of theorem~\ref{thm:perturbations1}. The well-posedness of the forward problem is proved by using the higher order fractional Poincar\'e inequality in theorem~\ref{thm:poincare} and interpolation inequality in non-homogeneous Sobolev spaces. In the case of coefficients with bounded fractional derivatives we also need the Kato--Ponce inequality in proving well-posedness. The higher order unique continuation of~$\fraclaplace$ in theorem~\ref{thm:ucplaplacian} together with well-posedness implies Runge approximation for equation~\eqref{eq:perturbations} (and for the adjoint equation of~\eqref{eq:perturbations}): one can approximate functions in $\widetilde{H}^s(\Omega)$ arbitrarily well by solutions of equation~\eqref{eq:perturbations}. We can prove uniqueness for the inverse problem by using the Runge approximation and suitable test functions in the Alessandrini identity which gives the relation between the DN maps~$\Lambda_{P_i}$ and the partial differential operators~$P_i$ in terms of exterior values~$f$ and solutions~$u_f$ of equation~\eqref{eq:perturbations} (and the adjoint equation of~\eqref{eq:perturbations}). It is important to notice that in theorems~\ref{thm:perturbations1} and~\ref{thm:perturbations2} we recover the coefficients~$a_\alpha$ uniquely and there is no gauge in contrast to the perturbed local Schr\"odinger equation~\cite{ISA-inverse-partial-differential, NSU95, SA-first-order-perturbations}.

Note that even though we consider partial differential operators in theorems~\ref{thm:perturbations1} and~\ref{thm:perturbations2}, the results apply for more general local linear operators. In fact, Peetre's theorem implies that any local linear operator $L\colon C_c^\infty(\Omega)\rightarrow C_c^\infty(\Omega)$ which satisfies $\spt(Lf)\subset\spt(f)$ for all $f\in C_c^\infty(\Omega)$ can be identified with a partial differential operator~\cite{MI-refinement-of-peetres-theorem, PEE-peetres-theorem}. Hence our results hold for any such local operator satisfying the assumptions in theorems~\ref{thm:perturbations1} and~\ref{thm:perturbations2}.

\section{Travel time tomography on Riemannian and Finsler manifolds: \cite{item:mixingray, item:randers}}
\label{sec:traveltimefinsler}

In the articles~\cite{item:mixingray, item:randers} we study the travel time tomography or boundary rigidity problem and its linearized versions on Riemannian manifolds and more general Finsler manifolds where the fiberwise inner product depends on direction. The basic idea of the boundary rigidity problem is illustrated in figure~\ref{fig:boundaryrigidity}. Suppose we have two Finsler norms~$F_1$ and~$F_2$ (which can be for example two Riemannian metrics) on a manifold~$M$ with boundary~$\partial M$. We assume that between any two boundary points $x, x'\in\partial M$ there is unique geodesic~$\gamma_i$ of the Finsler norm~$F_i$ going from~$x$ to~$x'$. The length of the geodesic~$\gamma_i$ with respect to~$F_i$ is denoted by~$L_{F_i}(\gamma_i)$ and it gives the (not necessarily symmetric) distance from~$x\in\partial M$ to~$x'\in\partial M$. The boundary rigidity problem is the following: if the Finsler norms~$F_1$ and~$F_2$ give the same distances between all boundary points $x, x'\in\partial M$, does it follow that $F_1=F_2$ up to a natural gauge?

\begin{figure}[htp]
\centering
\includegraphics[height=7.4cm]{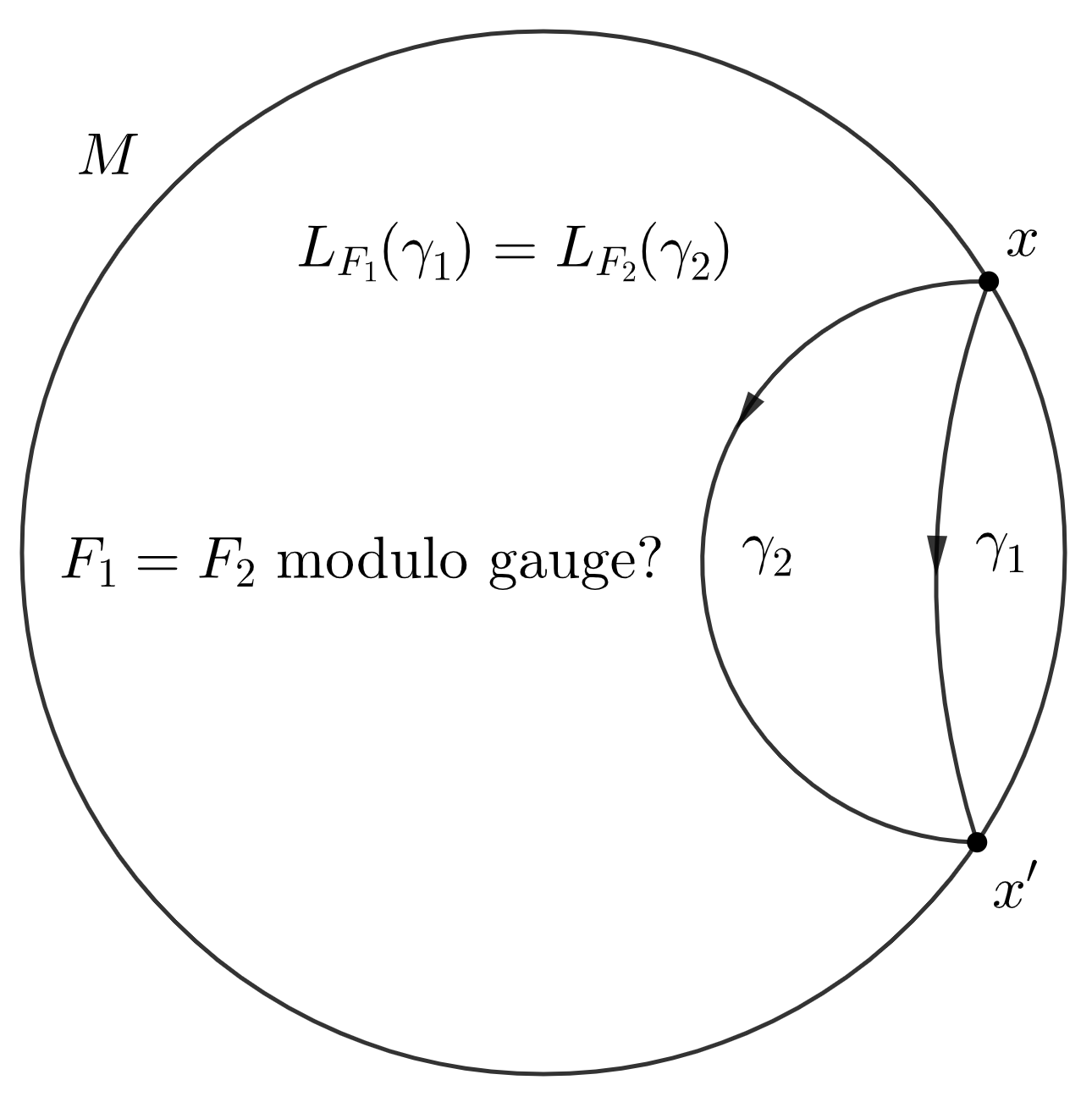}
\caption{An illustration of the boundary rigidity problem on Finsler manifolds $(M, F)$. Here $x, x'\in\partial M$ are two boundary points, $\gamma_1$ and $\gamma_2$ are the unique geodesics of the Finsler norms~$F_1$ and~$F_2$ connecting~$x$ to~$x'$, and $L_{F_i}(\gamma_i)$ denotes the length of the geodesic~$\gamma_i$ with respect to~$F_i$ (adapted from~\cite[Figure 1]{item:randers}).}
\label{fig:boundaryrigidity}
\end{figure}

If the Finsler norms~$F_i$ are induced by Riemannian metrics~$g_i$ (the fiberwise inner product does not depend on direction), then the natural gauge is a boundary preserving diffeomorphism: if $g_2=\Psi^*g_1$ where $\Psi\colon M\to M$ is a diffeomorphism such that $\Psi|_{\partial M}=\text{Id}$, then~$g_1$ and~$g_2$ give the same boundary distances. 
For a special class of non-reversible Finsler norms called Randers metrics the gauge is similar: if $F_1=F_g+\beta$ where~$F_g$ is a Finsler norm induced by the Riemannian metric~$g$ and $\beta$ is a 1-form whose norm with respect to~$g$ is small, then~$F_1$ and $F_2=\Psi^* F_1+\der\phi$ give the same boundary distances where $\Psi\colon M\to M$ is a diffeomorphism which is identity on the boundary and~$\phi$ is a scalar field vanishing on the boundary (and~$\der\phi$ is considered as a small perturbation to~$\Psi^*F_1$).

On Riemannian manifolds $(M, g)$ the linearization of the boundary rigidity problem leads to the geodesic ray transform~$\geod_m$ of symmetric (covariant) $m$-tensor fields~\cite{SHA-integral-geometry-tensor-fields}: if~$h$ and~$h'$ are two symmetric $m$-tensor fields such that $\geod_m h=\geod_m h'$, does it follow that $h=h'$ up to a natural gauge? When $m\geq 1$, the gauge is given by the derivative of a lower order tensor field: if $h'=h+\sigma\nabla v$ where~$h$ is symmetric $m$-tensor field, $v$ is an $m-1$-tensor field vanishing on the boundary (or at infinity) and~$\sigma\nabla$ is the symmetrized covariant derivative, then $\geod_m h=\geod_m h'$. Since the problem is linear ($\geod_m$ is a linear operator) it is enough to study the kernel of~$\geod_m$: if~$h$ is a symmetric $m$-tensor field such that $\geod_m h=0$, does it follow that $h=\sigma\nabla v$ where~$v$ is an $m-1$-tensor field vanishing on the boundary (or at infinity)? If this is true for all sufficiently regular symmetric $m$-tensor fields, we say that~$\geod_m$ is solenoidally injective (or s-injective).

In the article~\cite{item:mixingray} we study the mixed ray transform and more general mixing ray transforms on Riemannian manifolds. These integral transforms are generalizations of the geodesic ray transform and they arise in the linearization of the elastic travel time tomography problem~\cite{dESUZ-generic-uniqueness-mixed-ray, deSZ-mixed-ray, SHA-integral-geometry-tensor-fields}.
The main focus in the article~\cite{item:mixingray} is on the algebraic properties of mixing ray transforms and decompositions of tensor fields with respect to these transforms. We have various corollaries of a main idea how to study the kernel characterization and solenoidal injectivity of the mixing ray transforms using correct notion of symmetry and reduction. 

The mixing ray transform of $m$-tensor fields ($m\geq 1$) is defined as the composition $\geod_A h=(\geod_m\circ A)h$ where~$\geod_m$ is the geodesic ray transform of $m$-tensor fields and~$A$ is a smooth linear invertible map on $m$-tensor fields (see section~\ref{sec:notationtraveltime}). If $A$ is the identity map, then~$\geod_A$ reduces to the geodesic ray transform~$\geod_m$. One can think that the transform~$\geod_A$ first rotates the tensor field~$h$ and then takes the geodesic ray transform of the rotated $m$-tensor field~$Ah$.

The mixing ray transforms are matrix-weighted geodesic ray transforms and they have a different kind of kernel than the geodesic ray transform. We prove in the article~\cite{item:mixingray} that every $m$-tensor field~$h$ can be written as the direct sum $h=\widehat{\sigma}_A h+(h-\widehat{\sigma}_A h)$ where~$\widehat{\sigma}_A$ is the symmetrization map with respect to the transform~$\geod_A$ (see section~\ref{sec:notationtraveltime}) and $h-\widehat{\sigma}_A h\in\ker(\geod_A)$. Here~$\widehat{\sigma}_A h$ is the ``symmetric part" of~$h$ and $h-\widehat{\sigma}_A h$ is the ``trivial part" of~$h$ from the point of view of the transform~$\geod_A$. We show that if~$\geod_m$ is s-injective on symmetric $m$-tensor fields and $\geod_A h=0$, then $\widehat{\sigma}_A h=\widehat{\sigma}_A\nabla^A v$ for some $m-1$-tensor field~$v$ vanishing on the boundary (or at infinity) where $\nabla^A=A^{-1}\circ\nabla$ is the weighted covariant derivative associated to~$\geod_A$. This property is referred as the solenoidal injectivity of~$\geod_A$ and it allows us to write the kernel of~$\geod_A$ as the direct sum $\ker(\geod_A)=\im(\mathcal{H})\oplus\im(\widehat{\sigma}_A\nabla^A)$ where $\mathcal{H}=\id-\widehat{\sigma}_A$ is the projection onto the ``trivial part" of~$\ker(\geod_A)$ (see sections~\ref{sec:notationtraveltime} and~\ref{sec:mainresultstraveltime}).

In addition to solenoidal injectivity results we prove in the article~\cite{item:mixingray} numerous corollaries of the algebraic approach to mixing ray transforms and related transforms such as the mixed ray transform and the light ray transform. For example, we show that previous results for the light ray transform on Lorentzian manifolds and the mixed ray transform on simple Riemannian manifolds in~\cite{deSZ-mixed-ray, FIO-light-ray} can be seen as solenoidal injectivity results when we have a correct notion of symmetry. We also prove some stability results for the mixed ray transform, and show that the geodesic ray transform and the transverse ray transform together determine 1-forms uniquely on certain two-dimensional compact and non-compact manifolds.

In the article~\cite{item:randers} we study the boundary rigidity problem for certain non-reversible Finsler norms called Randers metrics. Finsler norms are non-negative functions on the tangent bundle $F\colon TM\rightarrow [0, \infty)$ so that for every $x\in M$ the map $y\mapsto F(x, y)$ is a positively homogeneous norm in~$T_xM$. The Finsler norm~$F$ is reversible, if $F(x, -y)=F(x, y)$ for all $x\in M$ and $y\in T_xM$. In this case the map $y\mapsto F(x, y)$ defines a norm in $T_x M$. In general, the distance function given by~$F$ is not necessarily symmetric in contrast to the Riemannian distance function. Finsler norms induce a fiberwise inner product which depends not only on position but also on direction. Riemannian metrics are a special case of reversible Finsler norms where the inner product does not depend on direction. 

Randers metrics are Finsler norms of the form $F=F_g+\beta$ where~$F_g$ is a Finsler norm induced by the Riemannian metric~$g$ and~$\beta$ is a 1-form whose norm with respect to~$g$ is small enough. Randers metrics are non-reversible since $F(x, -y)=F(x, y)$ for all $x\in M$ and $y\in T_x M$ if and only if $\beta\equiv 0$. Randers metrics arise naturally in Zermelo's navigation problem~\cite{BRS-zermelo-navigation, SHE-zermelo-riemannian}. Roughly saying, Zermelo's problem asks what is the shortest path in time for a moving object to travel from point~$A$ to point~$B$ when an external force field is acting on the object. Basic example is a ship which is sailing on a sea under the influence of wind or current.

In the article~\cite{item:randers} we prove two boundary rigidity results. If~$F$ is a Finsler norm and $x, x'\in\partial M$, denote by $d_F(x, x')$ the (non-symmetric) geodesic distance from~$x$ to~$x'$ (see section~\ref{sec:notationtraveltime}).  
The first theorem is the following: if~$F_1$ and~$F_2$ are Finsler norms of the form $F_i=F_{r, i}+\beta_i$ where~$F_{r, i}$ is a reversible Finsler norm and~$\beta_i$ is a closed 1-form ($\der\beta_i=0$) such that $d_{F_1}(x, x')=d_{F_2}(x, x')$ for all $x, x'\in\partial M$, then $\beta_2=\beta_1+\der\phi$ where~$\phi$ is a scalar field vanishing on the boundary and $d_{F_{r, 1}}(x, x')=d_{F_{r, 2}}(x, x')$ for all $x, x'\in\partial M$. This is done by using projective equivalence of the Finsler norms~$F_i$ and~$F_{r, i}$: since the 1-form~$\beta_i$ is closed~$F_i$ and~$F_{r, i}$ have the same geodesics as point sets, and the geodesics of~$F_i$ remain geodesics (as point sets) if their orientation is reversed. The second theorem is a boundary rigidity result for Randers metrics and it is a corollary of the first theorem: if~$F_1=F_{g_1}+\beta_1$ and~$F_2=F_{g_2}+\beta_2$ are Randers metrics where~$g_1$ and~$g_2$ are boundary rigid Riemannian metrics, $\beta_i$ is a closed 1-form and $d_{F_1}(x, x')=d_{F_2}(x, x')$ for all $x, x'\in\partial M$, then $F_2=\Psi^*F_1+\der\phi$ where~$\phi$ is a scalar field vanishing on the boundary and $\Psi\colon M\rightarrow M$ is a diffemorphism which is identity on the boundary. In other words, the equality of the boundary distances implies that the Randers metrics~$F_1$ and~$F_2$ are equal up to the natural gauge.
Using Zermelo's navigation problem we provide an application of the second theorem to seismology where the seismic wave propagates in a moving medium.

\subsection{Notation}
\label{sec:notationtraveltime}
Let us first go through the notation used in the article~\cite{item:mixingray}.
We follow the notation conventions of the references~\cite{LEE-smooth-manifolds, LEE-riemannian-manifolds, LRS-tensor-tomography-cartan-hadamard, PSU-tensor-tomography-progress, SHA-integral-geometry-tensor-fields}. We will use the Einstein summation convention so that every repeated index appearing both as a subscript
and superscript is implicitly summed over.

Let~$M$ be an $n$-dimensional smooth manifold where $n\geq 2$. We usually assume that~$M$ is compact and has a boundary~$\partial M$ or that~$M$ is non-compact without boundary.
If $(M, g)$ is a Riemannian manifold, we denote by~$K(x)$ the Gaussian curvature at $x\in M$. We say that a compact Riemannian manifold $(M, g)$ with boundary is simple (or that the Riemannian metric~$g$ is simple) if it is non-trapping (maximal geodesics have finite length), geodesics have no conjugate points and the boundary~$\partial M$ is strictly convex with respect to~$g$ (the second fundamental form on~$\partial M$ is positive definite). A compact simple manifold is always diffeomorphic to a ball. We say that a non-compact manifold $(M, g)$ without boundary is a Cartan--Hadamard manifold if it is simply connected, complete and its sectional curvature is nonpositive. Cartan--Hadamard manifolds are diffeomorphic to~$\R^n$ and basic examples are the Euclidean space and hyperbolic spaces.

Let $m\geq 1$. We denote by~$\X (T_mM)$ the space of all covariant $m$-tensor fields and $S_mM\subset \X(T_mM)$ is the space of symmetric covariant $m$-tensor fields. The notations~$C^\infty(T_m M):=C^\infty(\X(T_mM))$ and~$C^\infty(S_mM)$ mean that the corresponding tensor fields are smooth. The pointwise norm of a covariant $m$-tensor field~$h$ is $\abs{h}_{g_x}=\sqrt{g_x(h, h)}$ where~$g_x(\cdot, \cdot)$ is the fiberwise inner product of $m$-tensor fields.
We define the following sets of polynomially and exponentially decaying tensor fields which are mainly used on Cartan--Hadamard manifolds
\begin{align}\label{eq:cartanhadamardspaces}
E_{\eta}(T_mM)
&=
\{h\in C^1(T_mM):
\\&\qquad
\abs{h}_{g_x}\leq Ce^{-\eta d(x, o)} \ \text{for some} \ C>0\},
\\
E_{\eta}^1(T_mM)
&=
\{h\in C^1(T_mM):
\\&\qquad
\abs{h}_{g_x}+\abs{\nabla h }_{g_x}\leq Ce^{-\eta d(x, o)} \ \text{for some} \ C>0\},
\\
P_{\eta}(T_mM)&=\{h\in C^1(T_mM):
\\&\qquad
\abs{h}_{g_x}\leq C(1+d(x, o))^{-\eta} \ \text{for some} \ C>0\},
\\
P_{\eta}^1(T_mM)&=\{h\in C^1(T_mM):
\\&\qquad
\abs{h}_{g_x}\leq C(1+d(x, o))^{-\eta} \ \text{and}
\\&\qquad
\abs{\nabla h}_{g_x}\leq C(1+d(x, o))^{-\eta-1} \ \text{for some} \ C>0\}
\end{align}
where $o\in M$ is a fixed point and $\eta>0$.

Let $(M, g)$ be a non-trapping compact Riemannian manifold with boundary~$\partial M$. Let $x\in\partial M$ and $\xi\in T_xM$ be an inward-pointing unit vector. Denote by $\gamma_{x, \xi}$ the geodesic starting at~$x$ in the direction~$\xi$ and let $\tau(x, \xi)$ be the first time when the geodesic hits the boundary again.
The geodesic ray transform of a sufficiently regular $m$-tensor field~$h$ is defined as
\begin{equation}
\label{eq:grtbdry}
\geod_m h(x, \xi)=\int_0^{\tau(x, \xi)}h_{i_1\dotso i_m}(\gamma_{x, \xi}(t))\dot{\gamma}^{i_1}_{x, \xi}(t)\cdots\dot{\gamma}^{i_m}_{x, \xi}(t) \der t.
\end{equation}
Similarly, if $(M, g)$ is a Cartan--Hadamard manifold and $x\in M$ and $\xi\in T_xM$ has unit length, then we define
\begin{equation}
\geod_m h(x, \xi)=\int_{-\infty}^{\infty}h_{i_1\dotso i_m}(\gamma_{x, \xi}(t))\dot{\gamma}^{i_1}_{x, \xi}(t)\cdots\dot{\gamma}^{i_m}_{x, \xi}(t) \der t
\end{equation}
whenever the $m$-tensor field~$h$ decays rapidly enough at infinity. By completeness geodesics are defined on all times on Cartan--Hadamard manifolds.

We define $A\colon C^\infty(T_mM)\rightarrow C^\infty (T_mM)$ as a smooth linear invertible map on $m$-tensor fields which operates as
\begin{equation}\label{eq:mixingdefinition}
(Ah)_x(\xi_1,\dots,\xi_m) = h_x(A_1(x)\xi_1,\dots,A_m(x)\xi_m)
\end{equation}
where $\xi_i\in T_x M$ and each $A_i(x)$ is a linear bijection in~$T_xM$. Such map~$A$ is called a mixing of degree $m\geq 1$. If~$A$ is a mixing of degree~$m$, we define the mixing ray transform~$\geod_A$ by setting $\geod_A=\geod_m\circ A$ where~$\geod_m$ is the geodesic ray transform of $m$-tensor fields. On orientable two-dimensional manifolds an important special case of the mixing ray transforms is the mixed ray transform $\mixed=\geod_m\circ A_{k, l}$ where the components~$A_i$ of the mixing $A_{k, l}$ satisfy $A_i = \star$ when $i = 1,\dots, k$ and $A_i = \id$ when $i = k+1,\dots, k+l=m$. Here~$\star$ is the Hodge star operating on 1-forms (and hence on vector fields via the musical isomorphisms) and on orientable two-dimensional manifolds it corresponds to rotation by 90 degrees counterclockwise.

If~$A$ is a mixing of degree~$m$, we define the generalized symmetrization operator $\widehat{\sigma}_A=A^{-1}\circ\sigma\circ A$ where~$\sigma$ is the usual symmetrization of tensor fields. Then $\widehat{\sigma}_A$ is a projection onto~$A^{-1}(S_mM)$ and we have the direct decomposition $h=\widehat{\sigma}_A h+(h-\widehat{\sigma}_Ah)$ where $\widehat{\sigma}_Ah\in A^{-1}(S_mM)$ and $h-\widehat{\sigma}_Ah\in\ker(\geod_A)$. We denote by~$\nabla^A$ the weighted covariant derivative $\nabla^A=A^{-1}\circ\nabla$.
We say that the mixing ray transform~$\geod_A$ is s-injective on a compact Riemannian manifold $(M, g)$ with boundary, if for every $h\in C^\infty(T_mM)$ we have that $\geod_A h=0$ if and only if $\widehat{\sigma}_Ah=\widehat{\sigma}_A\nabla^A v$ for some $v\in C^\infty(S_{m-1}M)$ vanishing on the boundary.

Then we shortly introduce the additional notation used in the article~\cite{item:randers}; these basic notions of Finsler geometry can be found in~\cite{AL-global-aspects-finsler-geometry, BCS-introduction-finsler-geometry,  CS-riemann-finsler-geometry, SHE-lectures-on-finsler-geometry}.

Let $F\colon TM\rightarrow [0, \infty)$ be a Finsler norm and denote by~$F_r$ a reversible Finsler norm, i.e. $F_r(x, -y)=F_r(x, y)$ for all $x\in M$ and $y\in T_xM$. If~$g$ is a Riemannian metric, then it defines a reversible Finsler norm~$F_g$ as $F_g(x, y)=\sqrt{g_{ij}(x)y^iy^j}$. We denote by~$\beta$ a smooth 1-form and say that~$\beta$ is closed, if $\der\beta=0$ where~$\der$ is the exterior derivative of differential forms. We define the dual norm of~$\beta$ as $\aabs{\beta}_{F^*}=\sup_{x\in M}F^*(x, \beta_x)$ where~$F^*$ is the co-Finsler norm in~$T^*M$. More specifically, $F^*(x, \beta_x)=\sup_{y\in T_xM, F(x, y)=1}\beta_x(y)$. If~$F$ is a Finsler norm and~$\beta$ is a 1-form such that $\aabs{\beta}_{F^*}<1$, then $F+\beta$ also defines a Finsler norm. 

We define admissible Finsler norms as follows:~$F$ is admissible, if for any two points $x, x'\in\partial M$ there exists unique geodesic~$\gamma$ of~$F$ going from~$x$ to~$x'$ having finite length. 
When~$F$ is admissible, we define the map~$d_F(\cdot, \cdot)\colon\partial M\times\partial M\rightarrow [0, \infty)$ as $d_F(x, x')=L_F(\gamma)$ where~$L_F(\gamma)$ is the length of the geodesic~$\gamma$ with respect to~$F$. In general the map $d_F(\cdot, \cdot)$ is not symmetric.
We say that the Riemannian metrics $g_1$ and $g_2$ on~$M$ are boundary rigid, if $d_{g_1}(x, x')=d_{g_2}(x, x')$ for all $x, x'\in\partial M$ if and only if $g_2=\Psi^*g_1$ where $\Psi\colon M\rightarrow M$ is a diffeomorphism which is identity on the boundary.

\subsection{Main results}
\label{sec:mainresultstraveltime}
In the article~\cite{item:mixingray} we study linearized travel time tomography. We have numerous corollaries of the algebraic approach to mixing ray transforms and here we only present the most important results considering solenoidal injectivity. The first result says that s-injectivity of one mixing ray transform implies s-injectivity for all mixing ray transforms.

\begin{theorem}[{\cite[Corollary 3.4]{item:mixingray}}]
\label{thm:mixingraysinjectivity}
Let $m\geq 1$ and $(M, g)$ be a compact Riemannian manifold with boundary so that the transform~$I_A$ is s-injective for some~$A$ of degree~$m$. Then~$I_{\widetilde{A}}$ is s-injective for all~$\widetilde{A}$ of degree~$m$. 
\end{theorem}

Theorem~\ref{thm:mixingraysinjectivity} holds in all dimensions $n\geq 2$ and it is proved by using the definition of s-injectivity and the properties of mixings~$A$ and the corresponding projections~$\widehat{\sigma}_A$. The following theorem is a special case of theorem~\ref{thm:mixingraysinjectivity} in two dimensions.

\begin{theorem}[{\cite[Corollary 4.1]{item:mixingray}}]
\label{thm:sinjectivitycompact}
Let $m\geq 1$. Let $(M, g)$ be a compact two-dimensional orientable Riemannian manifold with boundary such that the geodesic ray transform is s-injective on $C^{\infty}(S_mM)$ and let $h\in C^{\infty}(T_mM)$.
Then $\mixed h=0$ if and only if $\widehat{\sigma}_{A_{k, l}}h=\widehat{\sigma}_{A_{k, l}}\nabla^{A_{k, l}} v$ for some $v\in C^{\infty}(S_{m-1}M)$ vanishing on the boundary~$\partial M$.
\end{theorem}

S-injectivity of the geodesic ray transform is known for example on compact simple surfaces~\cite{PSU-tensor-tomography-on-simple-surfaces} and on simply connected compact surfaces with strictly convex boundary and non-positive sectional curvature~\cite{PS-sharp-stability-nonpositive-curvature, SHA-integral-geometry-tensor-fields}.
Hence we obtain many new s-injectivity results for the mixed ray transform in two dimensions using theorem~\ref{thm:sinjectivitycompact}.
The assumption that $(M, g)$ is a two-dimensional orientable manifold is needed so that the mixed ray transform is well-defined, i.e. we can use the Hodge star~$\star$ to rotate vector fields.

The next theorem of the article~\cite{item:mixingray} shows that s-injectivity holds for the mixed ray transform also on certain non-compact Cartan--Hadamard manifolds.

\begin{theorem}[{\cite[Corollary 4.2]{item:mixingray}}]
\label{thm:c-hsinjectivity}
Let $(M, g)$ be a two-dimensional Cartan--Hadamard manifold and let $m\geq 1$.
The following claims are true:
\smallskip
\begin{enumerate}[(a)]
    \item\label{item:cartan1}
    Let $-K_0\leq K\leq 0$ for some $K_0>0$ and $h\in E^1_{\eta}(T_mM)$ for some $\eta>\frac{3}{2}\sqrt{K_0}$.
    Then $\mixed h=0$ if and only if $ \widehat{\sigma}_{A_{k, l}}h=\widehat{\sigma}_{A_{k, l}}\nabla^{A_{k, l}} v$ for some $v\in S_{m-1}M$ such that $v\in E_{\eta-\epsilon}(T_{m-1}M)$ for all $\epsilon>0$. 
    \medskip
    \item\label{item:cartan2}
    Let $K\in P_{\kappa}(M)$ for some $\kappa>2$ and $h\in P^1_{\eta}(T_mM)$ for some $\eta>2$.
    Then $\mixed h=0$ if and only if $ \widehat{\sigma}_{A_{k, l}}h=\widehat{\sigma}_{A_{k, l}}\nabla^{A_{k, l}} v$ for some $v\in S_{m-1}M\cap P_{\eta-1}(T_{m-1}M)$.
\end{enumerate}
\end{theorem}

Theorem~\ref{thm:c-hsinjectivity} follows from the corresponding s-injectivity result for the geodesic ray transform proved in~\cite{LRS-tensor-tomography-cartan-hadamard}. Before we can use the results in~\cite{LRS-tensor-tomography-cartan-hadamard} we show that the mixing~$A_{k, l}$ in the mixing ray transform $\mixed=\geod_m\circ A_{k, l}$ maps tensor fields in~$E^1_\eta(T_mM)$ to tensor fields in~$E^1_\eta (T_mM)$, and similarly tensor fields in~$P^1_\eta(T_mM)$ to tensor fields in~$P^1_\eta(T_mM)$. We do not need to assume orientability in theorem~\ref{thm:c-hsinjectivity} since Cartan--Hadamard manifolds are always orientable.

Theorem~\ref{thm:sinjectivitycompact} implies that we can write the kernel of the mixed ray transform on compact orientable surfaces with boundary admitting s-injectivity of the geodesic ray transform as the direct sum
\begin{equation}
\label{eq:decompositionkernel}
\ker(L_{k, l}|_{C^\infty(T_mM)})=\im(\mathcal{H}|_{C^\infty(T_mM)})\oplus \im(\widehat{\sigma}_{A_{k, l}}\nabla^{A_{k, l}}|_Y)
\end{equation}
where $\mathcal{H}=\id-\widehat{\sigma}_{A_{k, l}}$ is the projection onto the trivial part of $\ker(L_{k, l})$ and $Y=\{v\in C^\infty (S_{m-1}M):v|_{\partial M}=0\}$.
Similar decomposition as in~\eqref{eq:decompositionkernel} holds for non-compact Cartan--Hadamard manifolds by theorem~\ref{thm:c-hsinjectivity} using the sets of polynomially and exponentially decaying tensor fields.

In the article~\cite{item:randers} we study the non-linear travel time tomography or boundary rigidity problem. The following theorem is the first main result of the article~\cite{item:randers}.

\begin{theorem}[{\cite[Theorem 1.3]{item:randers}}]
\label{thm:rigidity1}
Let~$M$ be a compact and simply connected smooth manifold with boundary. For $i\in\{1, 2\}$ let $F_i=F_{r, i}+\beta_i$ be admissible Finsler norms where~$F_{r, i}$ is an admissible and reversible Finsler norm and~$\beta_i$ is a smooth closed 1-form such that $\aabs{\beta_i}_{F^*_{r, i}}<1$. Then the following are equivalent:
\begin{enumerate}[label=(\roman*)]
    \item\label{item:boundarydatathm1} $d_{F_1}(x, x')=d_{F_2}(x, x')$ for all $x, x'\in \partial M$.
    \medskip
    \item\label{item:gaugethm1} There is unique scalar field~$\phi$ vanishing on the boundary such that $\beta_2=\beta_1+\der\phi$, and $d_{F_{r, 1}}(x, x')=d_{F_{r, 2}}(x, x')$ for all $x, x'\in\partial M$.
\end{enumerate}
\end{theorem}
One can take~$F_r$ to be for example a simple Riemannian metric in theorem~\ref{thm:rigidity1} since they are admissible and reversible. 
Since~$F_r$ is reversible for any curve~$\gamma$ we can obtain~$L_{F_r}(\gamma)$ from the symmetric part and $\int_\gamma\beta$ from the antisymmetric part of the length functional~$L_F(\gamma)$. In other words, the data for~$\beta$ and~$F_r$ ``decouple". Closedness of the 1-form~$\beta$ is in essential role in proving theorem~\ref{thm:rigidity1}: $\der\beta=0$ implies that~$F_r$ and $F=F_r+\beta$ have the same geodesics up to orientation preserving reparametrizations, and geodesics of~$F$ remain geodesics as point sets when their parametrization is reversed. Simply connectedness of~$M$ implies that $\beta_i=\der\phi_i$ for some scalar field~$\phi_i$ and this fact also plays a role in the proof.

The next theorem is the second main result of the article~\cite{item:randers} and it gives a boundary rigidity result for certain Randers metrics.

\begin{theorem}[{\cite[Theorem 1.5]{item:randers}}]
\label{thm:rigidity2}
Let~$M$ be a compact and simply connected smooth manifold with boundary. For $i\in\{1, 2\}$ let $F_i=F_{g_i}+\beta_i$ be admissible Finsler norms where~$g_i$ is an admissible Riemannian metric and~$\beta_i$ is a smooth closed 1-form such that $\aabs{\beta_i}_{g_i}<1$. Assume that $(M, g_i)$ is boundary rigid. Then the following are equivalent:
\begin{enumerate}[label=(\alph*)]
    \item\label{item:boundarydatathm2} $d_{F_1}(x, x')=d_{F_2}(x, x')$ for all $x, x'\in \partial M$.
    \medskip
    \item\label{item:gaugethm2} There is unique scalar field~$\phi$ vanishing on the boundary and a diffeomorphism~$\Psi$ which is identity on the boundary such that $\beta_2=\beta_1+\der\phi$ and $g_2=\Psi^*g_1$.
    \medskip
    \item\label{item:gauge2thm2} There is unique scalar field~$\phi$ vanishing on the boundary and a diffeomorphism~$\Psi$ which is identity on the boundary such that $\beta_2=\Psi^*\beta_1+\der\phi$ and $g_2=\Psi^*g_1$.
\end{enumerate}
\end{theorem}

Theorem~\ref{thm:rigidity2} is proved by using theorem~\ref{thm:rigidity1} and the rigidity assumption on the Riemannian metrics~$g_i$. Theorem~\ref{thm:rigidity2} part~\ref{item:gauge2thm2} implies that $F_2=\Psi^*F_1+\der\phi$. Finsler norms satisfying such relation are sometimes called almost isometric Finsler norms and the diffeomorphism $\Psi\colon (M, F_2)\rightarrow (M, F_1)$ is called an almost isometry~\cite{CJ-functional-representation-almost-isometries,DJV-lipschitz-almost-isometries, HA-spacetimes-randers, JLP-almost-isometries-spacetimes}. We note that~$\Psi$ cannot be an isometry since this would require that $\Psi^*\beta_1=\beta_2$. Theorem~\ref{thm:rigidity2} can be seen as a generalization of the Riemannian boundary rigidity results to non-reversible Randers manifolds. If $n=2$, then one can take~$g_i$ to be a simple Riemannian metric in theorem~\ref{thm:rigidity2} since in two dimensions simple Riemannian metrics are boundary rigid~\cite{PU-simple-manifolds-boundary-rigidity}.

Theorem~\ref{thm:rigidity2} has the following application in seismology. Assume that $M=\overline{B}(0, R)\subset\R^n$ is a closed ball of radius $R>0$ equipped with the Riemannian metric $g=c^{-2}(r)e$ where $c=c(r)$ is a radial sound speed satisfying the Herglotz condition
\begin{equation}\label{eq:herglotz}
\frac{\der}{\der r}\bigg(\frac{r}{c(r)}\bigg)>0, \quad r\in [0, R],
\end{equation}
and~$e$ is the Euclidean metric. In addition, let us assume that~$g$ has no conjugate points so that~$g$ becomes a simple Riemannian metric~\cite{MO-herglotz-conjugate-points, SUVZ-travel-time-tomography}. Suppose that the seismic wave propagates in a moving medium which velocity field is given by the vector field~$W$. Using Zermelo's navigation problem and a first-order approximation we obtain that if the scaled flow field $W/c^2$ is irrotational ($\der (W/c^2)=0$), then one can uniquely determine the speed of sound~$c$ and the velocity field~$W$ up to potential fields from travel time measurements of seismic waves which are done on the boundary $\partial M=S^{n-1}(0, R)\subset\R^n$.

\bibliography{refs} 
\bibliographystyle{abbrv}

\end{document}